\documentclass[10pt]{article}

\newcommand{\deq}{:=} 
\usepackage{amsmath}
\usepackage{graphicx}
\usepackage{amsfonts}
\usepackage{amssymb}
\usepackage{subfig}
\usepackage{algorithm}
\usepackage{algorithmic}
\usepackage{mathtools,xspace}  
\mathtoolsset{showonlyrefs} 
\usepackage{float}
\usepackage{bbm}
\usepackage{epsfig}
\usepackage{framed}
\usepackage{caption}

\newcommand{\eg}[0]{\emph{e.g.,} }

\newcommand{\reals}{{\mathbb{R}}}

\def \expect{\mathbb{E}}
\DeclareMathOperator*{\argmin}{\arg\min}

\title{Regularized Non-Gaussian Image Denoising}
\date{}
\author{Albert K. Oh\footnotemark[2]
\and Rebecca M. Willett\footnotemark[2]}

\begin{document}
\maketitle

\renewcommand{\thefootnote}{\fnsymbol{footnote}}

\footnotetext[2]{University of Wisconsin-Madison, Department of Electrical and Computer Engineering \\ 1415 Engineering Drive, Madison, WI 53706.}

\renewcommand{\thefootnote}{\arabic{footnote}}

\begin{abstract}
In image denoising problems, one widely-adopted approach is to minimize a regularized data-fit objective function, where the data-fit term is derived from a physical image acquisition model. Typically the regularizer is selected with two goals in mind: (a) to accurately reflect image structure, such as smoothness or sparsity, and (b) to ensure that the resulting optimization problem is convex and can be solved efficiently. The space of such regularizers in Gaussian noise settings is well studied; however, non-Gaussian noise models lead to data-fit expressions for which entirely different families of regularizers may be effective. These regularizers have received less attention in the literature because they yield non-convex optimization problems in Gaussian noise settings. This paper describes such regularizers and a simple reparameterization approach that allows image reconstruction to be accomplished using efficient convex optimization tools. The proposed approach essentially modifies the objective function to facilitate taking advantage of tools such as proximal denoising routines. We present examples of imaging denoising under exponential family (Bernoulli and Poisson) and multiplicative noise.  \end{abstract}

\smallskip
\noindent \textbf{Keywords.} Image denoising, convex optimization, total variation, sparsity, proximal operators, image regularization

\pagestyle{myheadings}
\thispagestyle{plain}
\markboth{ALBERT OH AND REBECCA WILLETT}{REGULARIZED NON-GAUSSIAN IMAGE DENOISING}

\section{Introduction}
\label{sec:intro}
In image reconstruction, measured image data is always corrupted by noise. Active research addresses the problem of removing various types of noise from the observations based on models of both the data acquisition process and of underlying image structure. In many settings this noise is modeled using a Gaussian distribution because it is a good model of errors associated with common detectors like CCD arrays. However, many more specialized applications yield noise better modeled with alternative distributions, such as
Poisson \cite{snyder}, Gamma \cite{SARbook}, Rayleigh \cite{Sarti2005}, Bernoulli \cite{onebit}, and Geometric \cite{goyal2013}.
In these and other settings, a common denoising approach is to find an image which minimizes a regularized negative log-likelihood of the noisy observations, with regularizers adopted from the Gaussian image denoising literature. This paper shows that alternative families of regularizers, adapted to the noise model, can yield empirical performance gains. Furthermore, regularizers which are non-convex in Gaussian settings can be convex and quite effective in these non-Gaussian settings. 

More specifically, our interest is in recovering a true underlying image $x^* \in \mathbb{R}^p$ from the noise-corrupted observations $y \in \mathbb{R}^p$, where $y$ is drawn from a distribution parameterized by $x^*$ which models the noise, with the density denoted $p(y|x)$.
To obtain an estimate $\hat{x} \in \mathbb{R}^p$ of $x^*$, a common approach is to minimize a smooth, convex data-fit term $\phi(\cdot)$ regularized by a convex (and potentially non-smooth) regularization term $\rho(\cdot)$\footnote{Kernel and patch based methods, including bilateral filters, non-local means, and BM3D implicitly employ a {\em data-dependent} regularization function. Such approaches are beyond the scope of this paper.}. To formalize this notion, we consider
\begin{align}
\hat{x} = \argmin_x \ \{ \Phi(x,y) \deq \phi(x,y) + \tau \rho(x) \},
\end{align}
where $\tau>0$ is a tuning parameter.
In Gaussian noise settings, we often use $\phi(x,y) = \|x-y\|_2^2$, but more generally, a common choice for $\phi(x,y)$ is a term proportional to $-\log p(y|x)$, the negative log-likelihood of $y$ given $x$ according to a distribution selected to model noise statistics. In general, $\phi$ reflects the characteristics of the noise corrupting our observations.

The regularizer $\rho(\cdot)$ has been extensively studied and there are a few examples widely used in image recovery techniques. The choice of this regularizer depends on the assumptions made about the underlying image structure. Popular choices include the $\ell_1$ norm for coefficient sparsity in a wavelet basis or other dictionary, the total variation (TV) semi-norm for image gradient sparsity, and Huber-like functions which are akin to the $\ell_1$ norm but smooth.
In general, image processors choose a regularizer according to two desiderata:
\vspace{3mm}
\begin{itemize}
\item[(D1)] the objective function $\Phi(x,y)$ may be minimized efficiently; and 
\item[(D2)] the regularizer $\rho(x)$ accurately reflects image structure. 
\end{itemize}
\vspace{3mm} For instance, when the data-fit term is simply the sum of squared errors $\|x-y\|^2_2$, then one would choose a convex regularizer $\rho$.  {\bf In this paper, we describe several classes of image denoising problems where the noise is non-Gaussian and where {\em non-convex} regularization functions (which would be undesirable in Gaussian noise settings) not only admit efficient convex optimization algorithms, but also yield substantial empirical performance gains. }

As an illustrating example, consider two different denoising approaches, both based on total variation regularization \cite{osherTV}:
\begin{align}
\hat{x}_1 = \argmin_x \ -\log p(y|x) + \tau_1 \| x \|_{\rm TV}
\label{eq:theta1}
\end{align}
and
\begin{align}
\hat{x}_2 = \argmin_x \ -\log p(y|x) + \tau_2 \| \log x \|_{\rm TV},
\label{eq:theta2}
\end{align}
where
\begin{align}
\| x \|_\mathrm{TV} := & \sum_{i=1}^{p_1-1} \sum_{j=1}^{p_2-1} \sqrt{(x_{i,j} - x_{i+1,j})^2 + (x_{i,j} - x_{i,j+1})^2}  \\
& + \sum_{i=1}^{p_1-1} |x_{i,p_2} - x_{i+1,p_2}| + \sum_{j=1}^{p_2-1} |x_{p_1,j} - x_{p_1,j+1}|
\label{eq:TV}
\end{align}
is known as the discrete isotropic TV \cite{FGP}, where in a slight abuse of notation we consider $x$ to be an $p_1 \times p_2$ image matrix.

Comparing the two optimizations, the main difference is that the TV norm is applied to either the image \cite{SPIRAL} or the $\log$ of the image (assuming all the pixel values are positive) \cite{logTV}. If we think about the TV norm as promoting sparse image gradients, then it is clear that applying the TV norm to the log-intensity image is a sensible reflection of image structure, and may allow for larger non-zero gradient magnitudes than the TV norm applied directly to the intensity. If the noise were Gaussian, so that $-\log p(y|x) \propto \|y-x\|_2^2$, the optimization problem in \eqref{eq:theta2} would be non-convex, making this regularizer undesirable. However, in the Poisson noise setting, where $-\log p(y|x) \propto \sum_{i=1}^p x_i + y_i \log x_i$, a straightforward reparameterization (as described in \cite{logTV}) leads to a convex optimization problem.

{\em In other words, the shape of the log-likelihood function helps determine the set of regularization functions that satisfy (D1) and (D2). This means that the best choice of regularizer for a given problem may depend on the noise model and hence may be quite different from what one would choose in a Gaussian noise setting.}

\subsection{Organization of the paper}
Section~\ref{sec:form} contains an explicit formulation of the the problem, and Section~\ref{sec:approach} details our approach. Section \ref{sec:compstat} presents how reparameterization can affect computational-statistical tradeoffs.
In Section \ref{sec:expfam}, we introduce a generalization for image reparameterization for the exponential family of distributions; in Sections \ref{sec:bernoulli} and \ref{sec:poisson} we consider two specific examples of Bernoulli and Poisson noise, respectively. In Section \ref{sec:gamma}, we consider an example of reparameterization for multiplicative (speckle) noise. In Section \ref{sec:experiments}, we present experimental results in selecting the best reparameterization (in terms of RMSE) over a bank of images.
We conclude the paper in Section \ref{sec:conclusion}.

\section{Problem formulation}
\label{sec:form}

In choosing a good regularizer, one might search over a space of candidate regularizers such that the objective function $\Phi(x,y) = \phi(x,y) + \rho(x)$ is convex and choose a regularizer $\rho$ from this space which yields accurate estimates across a collection of images.
However, this space can be very difficult to specify and overly restrictive in practice, especially since $\phi(x,y)$ is not necessarily convex in $x$. We consider a tractable alternative in this paper. Specifically, we start with a {\em fixed}, {\em known} family of convex image regularizers $\rho$, denoted $\cal{R}$. For each $\rho \in \cal{R}$, we search over a space of candidate reparameterization functions $f: \reals^p \mapsto \reals^p$ such that the objective function
$$\Phi(f(\theta),y) = \phi(f(\theta),y) + \tau\rho(\theta)$$ is convex in $\theta$ for all $\tau>0$.
For any such $f$, given a new parameter $\theta \in \reals^p$, $x=f(\theta) \in \reals^p$ is a reparameterized version of the image intensity. 
In other words, instead of searching over regularizers that satisfy (D1) and (D2), we search over corresponding reparameterizations of the log likelihood. This approach guarantees that we will be able to leverage off-the-shelf convex optimization solvers that rely on proximal operators for a regularization function $\rho$ \cite{fista, sparsa}.

More specifically, let $\cal{R}$ be a collection of convex regularization functions $\rho:\reals^p\mapsto\reals$ that reflect desirable image structure (\eg sparsity, smoothness, etc.). Let $\Theta \subseteq \reals^p$ be a feasible parameter space. Define the set of reparameterization functions 
\begin{equation}
\label{eq:F}
\mathcal{F} = \{f:\Theta\mapsto\reals^p: \phi(f(\cdot),y) \mathrm{\ convex \ and \ smooth } \}.
\end{equation}
Given some $f \in \cal{F}$, $\rho \in \cal{R}$, and $\tau > 0$, we may estimate our image via:
\begin{align} 
\label{eq:obj_theta}
\hat{\theta}_{f,\rho,\tau} &= \argmin_\theta \ \{ \Phi(f(\theta),y) = \phi(f(\theta),y) + \tau \rho(\theta) \} \\
\hat{x}_{f,\rho,\tau} &= f(\hat{\theta}_{f,\rho,\tau}),
\label{eq:estimator}
\end{align}
where it is implied that $\hat{\theta}_{f,\rho,\tau}$ and $\hat{x}_{f,\rho,\tau}$ are functions of the noise $y$, but we omit writing this explicitly throughout this paper for brevity.

Our goal is to choose a regularizer $\rho \in \cal{R}$ and reparameterization $f \in \cal{F}$ according to
$$(f^*,\rho^*) = \argmin_{(f,\rho) \in \cal{F} \times \cal{R}} R(f,\rho)$$
where $R$ measures the {\em risk}
\begin{align}
R(f,\rho) \deq \expect_{X,Y} \left[ \min_{\tau > 0} \|\hat{x}_{f,\rho,\tau}-X\|_2 / \|X\|_2 \right],
\label{eq:truerisk}
\end{align}
where the expectation is with respect to both the noise in the data $y$ (of which $\hat{x}_{f,\rho,\tau}$ is an implicit function) and a distribution over the space of images.

Notice that because the variable we are optimizing over is $\theta$, we are directly regularizing the new parameter $\theta$ and not the image. The interpretation is that by changing the reparameterization function $f$, we change the relationship between the image $x$ and the parameter vector $\theta$, and thus implicitly change how $x$ is regularized. To see this directly, note that if $f$ is invertible, \eqref{eq:estimator} is equivalent to
\begin{align}
\hat{x}_{f,\rho,\tau} = \argmin_x \ \{ \Phi(x,y) = \phi(x,y) + \tau \rho(f^{-1}(x)) \}.
\label{eq:implicit_reg}
\end{align}
We note that it is not required for $f$ to be invertible in general, but it is useful in establishing this interpretation. From this perspective, it is clear that the regularizer applied to the image is $\rho(f^{-1}(x))$, and so choosing a reparameterization function $f$ implicitly selects an associated regularizer. Furthermore, the implicit regularization function $\rho(f^{-1}(\cdot))$ is not necessarily convex, but convex optimization methods can still be used to denoise the image.

{\em The key point is that the form of the log-likelihood will govern the contents of $\cal{F}$ and the specific risk-minimizing reparameterization and associated regularizer.}
Conventional wisdom tells us that a ``good'' regularizer should simply be convex and reflect desirable image structure, and once a good regularizer is found it can be used with any noise model. In contrast, we will see examples where regularizers considered good according to this conventional wisdom are suboptimal. Furthermore, regularizers identified in our framework will often be non-convex functions of the image (and hence undesirable in Gaussian noise settings) but still leading to convex optimization formulations for non-Gaussian image estimation.

\section{Proposed approach}
\label{sec:approach}

Our goal is to minimize the risk defined in \eqref{eq:truerisk}. Since in general we do not have access to the distribution over the space of images, we minimize an empirical estimate of risk using a collection of $n$ representative images, denoted $x^{*(1)}, \ldots, x^{*(n)}$. For the $i^{th}$ image, we generate $m$ noisy realizations $y^{(i,1)}, \ldots, y^{(i,m)}$, whose reconstructions are denoted
$\hat{x}_{f,\rho,\tau}^{(i,1)}, \ldots, \hat{x}_{f,\rho,\tau}^{(i,m)}$, via \eqref{eq:estimator}. Depending on the noise distribution assumptions, we select an appropriate data-fit term $\phi$ based on the negative log-likelihood. Note that often the negative log-likelihood itself will be convex; however, alternative formulations of $\phi$ may be desired when it is non-convex; we present an example of this in Section~\ref{sec:gamma}. For a particular $(f,\rho)$ pair, we then compute the empirical risk
\begin{align}
\hat{R}(f,\rho) = \frac{1}{nm} \sum_{i=1}^n \min_{\tau>0} \sum_{j=1}^m \|\hat{x}_{f,\rho,\tau}^{(i,j)}-x^{*(i)}\|_2 / \|x^{*(i)} \|_2,
\label{eq:relative_risk}
\end{align}
which is essentially a metric that averages the normalized RMSE of $n$ image examples with $m$ noise realizations each, noting that $\tau$ is chosen to minimize the average (across noise realizations) RMSE for each image.
Finally, we are then interested in choosing a reparameterization $\hat{f}$ that minimizes this empirical risk, namely
\begin{align}
\hat{f} = \underset{f \in \mathcal{F}}{\argmin} \ \hat{R}(f).
\end{align}

In order to compute the estimates $\hat{x}_{f,\rho,\tau}^{(i,j)}$, we adopt a proximal gradient method \cite{sparsa}, where a key step in iteration $t$ is to solve the proximal denoising problem 
\begin{align}
\label{eq:prox}
\theta^{t+1}_{f,\rho,\tau} = \argmin_\theta \frac{1}{2} \|\theta - s_t\|^2_2 + \frac{\tau}{\alpha_t} \rho(\theta)
\end{align}
for appropriately specified gradient descent step $s_t$ and (potentially adaptive) step size $\alpha_t$. Because our strategy is to search over reparameterization functions, off-the-shelf proximal gradient methods can be used to compute $\hat{x}_{f,\rho,\tau}$ quickly and accurately \cite{fista, wajs2006, sparsa, rapid}.

\section{Computational-statistical tradeoffs}
\label{sec:compstat}

The proposed methodology can be used not only to improve empirical RMSE performance, but also to exploit
computational-statistical tradeoffs \cite{Bruer2014,chandrasekaran2013,shender2013}. Take the example of Poisson denoising in Section \ref{sec:intro}, comparing \eqref{eq:theta1} and \eqref{eq:theta2}. In light of Section \ref{sec:form}, we can rewrite \eqref{eq:theta2} equivalently as
\begin{align} 
\label{eq:obj_theta_ex}
\hat{\theta}_2 &= \argmin_\theta \ -\log p(y|\exp(\theta)) + {\tau}_2' \| \theta \|_\mathrm{TV} \\
\hat{x}_2 &= \exp(\hat{\theta}_2),
\label{eq:estimator_ex}
\end{align}
where $\exp(\cdot)$ is applied element-wise, which is how we would perform the ``logarithmic regularization" of \eqref{eq:theta2} in practice: first reparameterize the image with an exponential transform to obtain the estimate $\hat{\theta}_2$, then obtain the image estimate through the reparameterization. With this optimization problem, the data-fit term $-\log p(y|\exp \theta) \propto \sum_{i=1}^p \exp \theta_i - y_i\theta_i$ grows without bound (i.e. does not have a Lipschitz gradient) as $\theta \rightarrow \infty$ but is otherwise well-behaved (i.e. does have a Lipschitz gradient). On the other hand, $-\log p(y|x) \propto \sum_{i=1}^p x_i - y_i \log x_i$ grows without bound as $x \rightarrow 0$ but is otherwise well-behaved. So we see that in both cases, we have data-fit terms whose gradients are not globally Lipschitz, but are locally Lipschitz on some domain.

A Lipschitz gradient is important in characterizing the convergence rate of many convex optimization methods
\cite{fista}. In general, we can say that an algorithm applied to a function with a smaller Lipschitz constant on its gradient will converge much faster than one with a larger Lipschitz constant; as a general intuition, having a larger Lipschitz constant (i.e. greater curvature) requires smaller steps in a gradient descent setting, leading to slower convergence. This leads to a practical problem: many of the data-fit terms resulting from commonly-used noise distributions, and not just Poisson, often have regions with unbounded curvature.
The class of reparameterizations $\cal{F}$ described above can be constrained to yield data-fit terms with globally Lipschitz gradients, which in turn increases convergence speed; in this case $f(\theta)$ is a function that is selected to make regions of $\phi(f(\theta),y)$ with unbounded curvature in $f(\theta)$ have bounded curvature in $\theta$.

In settings where convergence speed is critical, we can use this idea to modify the objective function used to select a reparameterization and reguarization pair. Specifically, let $L(\phi(f,y))$ denote the Lipschitz constant of the gradient of $\phi(f(\cdot),y)$. Then instead of minimizing $R(f,\rho)$ as in \eqref{eq:truerisk}, we might compute
$$(f^*,\rho^*) = \argmin_{(f,\rho) \in \cal{F} \times \cal{R}} R(f,\rho) + \lambda (L(\phi(f,y)))^{-1},$$
where $\lambda > 0$ controls the computational-statistical tradeoff.
This idea will be demonstrated in the results presented in Section \ref{sec:experiments}.

\section{Reparameterization for the exponential family}
\label{sec:expfam}

Noise from the exponential family is commonly seen in imaging problems, so it is useful to present a general $\phi(f(\theta),y)$ function for this class of problems. Distributions in the exponential family can be written in the form
\begin{align}
p(y | \eta) = \prod_{i=1}^p h(y_i) \exp\left[ \eta_i T(y_i) - A(\eta_i) \right],
\end{align}
where $h(\cdot)$ is the base measure (a constant to scale the distribution so that it integrates to 1), $\eta \in \mathbb{R}^p$ are the natural parameters that depend on distribution parameters, $T(\cdot)$ is the sufficient statistic, and $A(\cdot)$ is the log-partition function \cite{ghosh2005}. The negative log-likelihood of this family is then proportional to
\begin{align}
\tilde{\phi}(\eta,y) = \sum_{i=1}^p A(\eta_i) - \eta_i T(y_i).
\label{eq:phi_eta}
\end{align}
Equivalently, we can rewrite the data-fit term \eqref{eq:phi_eta} in terms of the image intensities $x$, whose relationship to the natural parameter is $\eta = g(x)$,  where $g: \reals^p \mapsto \reals^p$ is an invertible function that maps from the image $x$ to the natural parameter $\eta$, defined for all distributions in the exponential family; we note that $g$ cannot always be applied element-wise, but in the specific examples to follow, we use distributions whose $g$ functions can be applied element-wise. We may then write the negative log-likelihood as
\begin{align}
\phi(x,y) = \sum_{i=1}^p A([g(x)]_i) - [g(x)]_i T(y_i).
\label{eq:phi_mu}
\end{align}

To develop our family of reparameterizations $\cal{F}$, we consider 
functions which bridge between the data fit terms in \eqref{eq:phi_eta} and \eqref{eq:phi_mu} via
\begin{align}
\phi(f(\theta),y) = \sum_{i=1}^p A([g(f(\theta))]_i) - [g(f(\theta))]_i T(y_i).
\label{eq:general_datafit}
\end{align}

Our task now is to form a class of reparameterization functions $\cal{F}$ where the following criteria are met:
\begin{enumerate}
\item for each $f \in \cal{F}$, the corresponding $\phi(f(\cdot),y)$ is continuous, convex, and has Lipschitz gradients,
\item regularizing in the new parameter space yields strong empirical performance; and
\item we may efficiently explore the space $\cal{F}$.
\end{enumerate}
To accomplish this, given $a, b \in \reals$, define
\begin{align}
f_{a,b}(\theta) = g^{-1}(a\theta + b).
\label{eq:f_ab}
\end{align}
We  define a class $\mathcal{F}$ where $[f(\theta)]_i$ is either $[f_{a,b}(\theta)]_i$ or $\theta_i$ depending on which has a smaller partial derivative with respect to $\theta_i$. Formally,
\begin{equation}
\begin{aligned}
\mathcal{F} = \left\{ \vphantom{\begin{cases} a \\ a \end{cases}} \right. f \ \text{continuous over}& \ \Theta: \exists a,b \in \reals \ \text{s.t.} \\
&[f(\theta)]_i = \begin{cases} [f_{a,b}(\theta)]_i, & \frac{\partial [f_{a,b}(\theta)]_i}{\partial \theta_i} < 1 \\
\theta_i, & \text{otherwise} \end{cases}, \ i = 1,\cdots,p \left. \vphantom{\begin{cases} a \\ a \end{cases}} \right \}.
\label{eq:family}
\end{aligned}
\end{equation}
Note that for $f \in \mathcal{F}$,
\begin{equation}
\begin{aligned}
[g(f(\theta))]_i = \begin{cases} a\theta_i+b, & \frac{\partial [f_{a,b}(\theta)]_i}{\partial \theta_i} < 1 \\
[g(\theta)]_i, & \text{otherwise} \end{cases}.
\label{eq:general_reparam}
\end{aligned}
\end{equation}

We note that this $f(\theta)$ reparameterization does not always yield a convex data-fit term for all distributions in the exponential family, though in what follows, we explore examples that are convex. We also note that the family defined in \eqref{eq:family} is only one possible family and others are certainly possible, though we found this family to be useful for the two examples in the following subsections. In both of the examples to follow, it is straightforward to show that the reparameterized data-fit terms we propose are convex and have Lipschitz gradients through basic algebra (convexity by showing that the diagonal Hessian matrix has diagonal elements $\frac{\partial^2 \phi(f(\theta),y)}{\partial \theta_i^2}$ that are nonnegative, and Lipschitz gradients by showing $\frac{\partial^2 \phi(f(\theta),y)}{\partial \theta_i^2}$'s are finite); continuity is guaranteed by construction.

Table \ref{tab:exp_fam} lists the relevant functions for the Bernoulli and Poisson distributions, which will be discussed in the following.

\begin{table}[ht]
\begin{center}
\begin{tabular}{ c | c c c c}
 & $A(\cdot)$ & $T(y)$ & $g(\cdot)$ & $g^{-1}(\cdot)$ \\ \hline
Bernoulli & $\log(1+\exp(\cdot))$ & $y$ & $\log(\frac{\cdot}{1-\cdot})$ & $\frac{\exp(\cdot)}{1+\exp(\cdot)}$ \\
Poisson & $\exp(\cdot)$ & $y$ & $\log(\cdot)$ & $\exp(\cdot)$ \\
\end{tabular}
\caption{Table of relevant functions for Bernoulli and Poisson distributions (log-partition $A$, sufficient statistic $T$, natural parameter mapping $g$, and inverse mapping $g^{-1}$). In these specific examples chosen for this paper, all the operations in this table can be applied element-wise.}
\label{tab:exp_fam}
\end{center}
\end{table}

\subsection{Bernoulli noise (noisy one-bit sensors)}
\label{sec:bernoulli}

One example of a denoising problem that fits into this framework is in quantization noise removal. Consider an example of an unreliable sensor system where measurements are known to be imprecise, and so the observations are taken to be binary \cite{onebit}. In this one-bit quantization noise case, the sensor can be modeled with a dithering and quantization so that
\begin{equation}
\begin{aligned}
y_i = I_{x_i > u_i} \deq \begin{cases} 1, & \mathrm{if} \ x_i > u_i \\
0, & \mathrm{otherwise} \end{cases},
\end{aligned}
\end{equation}
where $x_i \in [0,1]$ and corresponds to relative pixel intensities, and $u_i$ is a threshold drawn from Uniform(0,1). This can be modeled with a Bernoulli distribution, namely $y \sim \mathrm{Bernoulli}(x^*)$, which falls into the exponential family.

Adopting the strategy above in \eqref{eq:general_datafit} and applying Table \ref{tab:exp_fam}, the Bernoulli distribution results in the data-fit term
\begin{align}
\phi(f(\theta),y) = \sum_{i=1}^p \log(1+\exp[g(f(\theta))]_i) - y_i [g(f(\theta))]_i.
\label{eq:datafit_bernoulli}
\end{align}
Before defining $f(\theta)$ for the Bernoulli case, we note that we use different reparameterizations depending on if $\theta_i \le 0.5$ or $\theta_i > 0.5$; this is a small extension of \eqref{eq:family}. For this particular case we specify a family of reparameterizations as:
\begin{equation}
\begin{aligned}
{\mathcal{F}} = \left \{ \vphantom{\begin{cases} a \\ a \\ a \end{cases}} \right. f \ \text{continuous} & \ \text{over} \ \Theta: \exists a,b_1,b_2 \in \reals \ \text{s.t.} \\ 
&[f(\theta)]_i = \begin{cases} [f_{a,b_1}(\theta)]_i & \frac{\partial [f_{a,b_1}(\theta)]_i}{\partial \theta_i} < 1 \ \text{and} \ \theta_i \le 0.5 \\
 [f_{a,b_2}(\theta)]_i & \frac{\partial [f_{a,b_2}(\theta)]_i}{\partial \theta_i} < 1 \ \text{and} \ \theta_i > 0.5 \\
\theta_i & \text{otherwise} \end{cases}, \ i = 1,\cdots,p \left. \vphantom{\begin{cases} a \\ a \\ a \end{cases}} \right \}.
\label{eq:family_twosided}
\end{aligned}
\end{equation}
Note that for any function in this class, $g(f(\theta))$ has the form:
\begin{equation}
\begin{aligned}
[g(f(\theta))]_i = \begin{cases} a\theta_i+b_1, & \frac{\partial [f_{a,b_1}(\theta)]_i}{\partial \theta_i} < 1 \ \text{and} \ \theta_i \le 0.5 \\
a\theta_i+b_2, & \frac{\partial [f_{a,b_2}(\theta)]_i}{\partial \theta_i} < 1 \ \text{and} \ \theta_i > 0.5 \\
[g(\theta)]_i, & \text{otherwise} \end{cases}.
\label{eq:general_reparam_twosided}
\end{aligned}
\end{equation}
In the above, there are three degrees of freedom in $a$, $b_1$, and $b_2$. However, there is effectively only one because for any triple $(a,b_1,b_2)$ satisfying the continuity constraints, there exists some $k \in (0,0.5)$ such that
\begin{equation}
\begin{aligned}
[g(f(\theta))]_i = \begin{cases} \frac{1}{k(1-k)}\theta_i+\log(\frac{k}{1-k})-\frac{1}{1-k}, & \theta_i \le k\\
\log(\frac{\theta_i}{1-\theta_i}), & k < \theta_i \le 1-k \\
\frac{1}{k(1-k)} \theta_i+\log(\frac{1-k}{k}) - \frac{1}{k}, &\theta_i > 1-k \end{cases},
\end{aligned}
\end{equation}
where image intensities are remapped via
\begin{equation}
\begin{aligned}
[f(\theta)]_i \deq \begin{cases} \left[ 1 + \frac{1-k}{k} \exp \left( \frac{-1}{k(1-k)} \theta_i + \frac{1}{1-k} \right) \right]^{-1}, & \theta_i \le k\\
\theta_i, & k < \theta_i \le 1-k \\
\left[ 1 + \frac{k}{1-k} \exp \left( \frac{-1}{k(1-k)} \theta_i + \frac{1}{k} \right) \right]^{-1}, &\theta_i > 1-k \end{cases},
\label{eq:bern_f}
\end{aligned}
\end{equation}
so that the family of reparameterizations can be formally rewritten as 
\begin{equation}
\begin{aligned}
{\mathcal{F}} = \left \{ \vphantom{\begin{cases} \left[ 1 + \frac{1-c}{c} \exp \left( \frac{-1}{c(1-c)} \theta_i + \frac{1}{1-c} \right) \right]^{-1} \\ a \\ \left[ 1 + \frac{1-c}{c} \exp \left( \frac{-1}{c(1-c)} \theta_i + \frac{1}{1-c} \right) \right]^{-1} \end{cases}} \right. 
f: \Theta &\rightarrow \reals^p: \exists k \in (0,0.5) \ \text{s.t.} \\ 
&[f(\theta)]_i \deq \begin{cases} \left[ 1 + \frac{1-k}{k} \exp \left( \frac{-1}{k(1-k)} \theta_i + \frac{1}{1-k} \right) \right]^{-1}, & \theta_i \le k\\
\theta_i, & k < \theta_i \le 1-k \\
\left[ 1 + \frac{k}{1-k} \exp \left( \frac{-1}{k(1-k)} \theta_i + \frac{1}{k} \right) \right]^{-1}, &\theta_i > 1-k \end{cases}, \ i = 1,\cdots,p \left. \vphantom{\begin{cases} \left[ 1 + \frac{1-c}{c} \exp \left( \frac{-1}{c(1-c)} \theta_i + \frac{1}{1-c} \right) \right]^{-1} \\ a \\ \left[ 1 + \frac{1-c}{c} \exp \left( \frac{-1}{c(1-c)} \theta_i + \frac{1}{1-c} \right) \right]^{-1} \end{cases}} \right \}.
\end{aligned}
\end{equation}

\begin{framed}
Summarizing, in the Bernoulli case we propose the data-fit expression
\begin{align}
\phi(f(\theta),y) = \sum_{i=1}^p \phi_i(f(\theta),y),
\label{eq:bernoulli_final}
\end{align}
where
\begin{align}
\phi_i(f(\theta),y) \deq \begin{cases} \log \left(1 + \frac{k}{1-k} \exp(\frac{\theta_i-k}{k(1-k)}) \right) - y_i \left[ \frac{\theta_i-k}{k(1-k)}+\log(\frac{k}{1-k}) \right], & \theta_i \le k\\
\log(\frac{1}{1-\theta_i}) - y_i \log(\frac{\theta_i}{1-\theta_i}), & k < \theta_i \le 1-k \\
\log \left(1 + \frac{1-k}{k} \exp(\frac{\theta_i+k-1}{k(1-k)}) \right) - y_i \left[ \frac{\theta_i+k-1}{k(1-k)}+\log(\frac{1-k}{k}) \right], &\theta_i > 1-k \end{cases}
\end{align}
for some $k \in (0,0.5)$.
\end{framed}

It is straightforward to show that the reparameterized data-fit term $\phi(f(\theta),y)$ we propose is convex and has Lipschitz gradients through algebra (convexity by showing that the diagonal Hessian matrix has diagonal elements $\frac{\partial^2 \phi(f(\theta),y)}{\partial \theta_i^2}$ that are nonnegative, and Lipschitz gradients by showing $\frac{\partial^2 \phi(f(\theta),y)}{\partial \theta_i^2}$'s are finite).

Our proposed $f(\theta)$ is an invertible function. Invertibility is not required in our setting, but it is useful in understanding the implicit regularization function $\rho(f^{-1}(x))$ seen in \eqref{eq:implicit_reg}. The inverse function is
\begin{align}
[f^{-1}(x)]_i = \begin{cases} \left[ \frac{1}{1-k} - \log\left(\frac{k}{1-k}(\frac{1}{x_i}-1)\right) \right]k(1-k), & x_i \le k\\
x_i, & k < x_i \le 1-k \\
\left[\frac{1}{k} - \log\left(\frac{1-k}{k}(\frac{1}{x_i}-1)\right)\right]k(1-k), &x_i > 1-k \end{cases}.
\label{eq:bern_implicit}
\end{align}
For an intuition behind the regularization, Figure \ref{fig:reg_bern} illustrates how the implicit regularizer $\| f^{-1}(x) \|_\mathrm{TV}$ behaves.

\newlength{\figwidth}

\setlength{\figwidth}{0.33\textwidth}
\begin{figure}[t]
\begin{center}
\subfloat[Intensity reparam.]{\includegraphics[width=\figwidth]{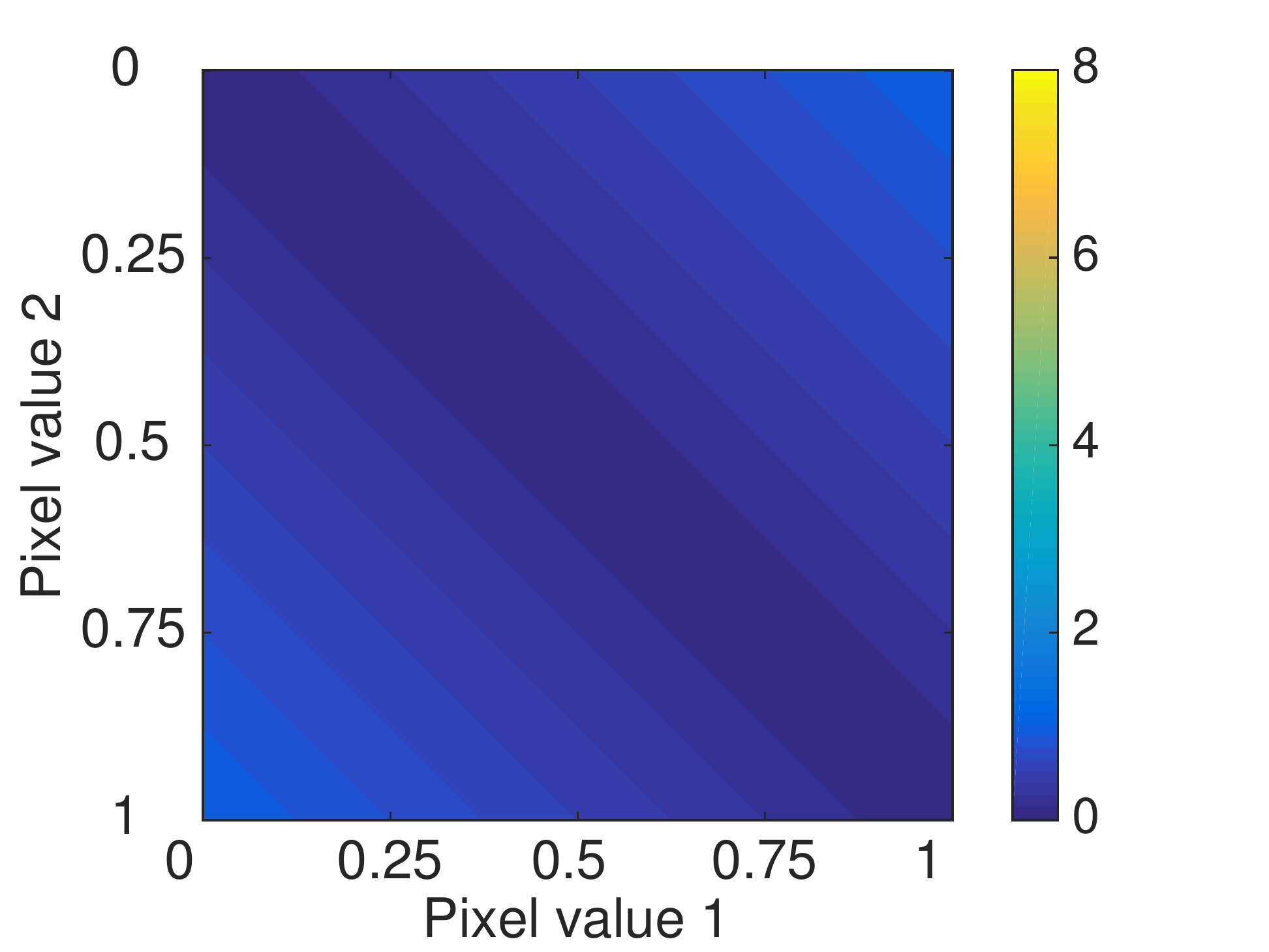}} ~
\subfloat[Natural reparam.]{\includegraphics[width=\figwidth]{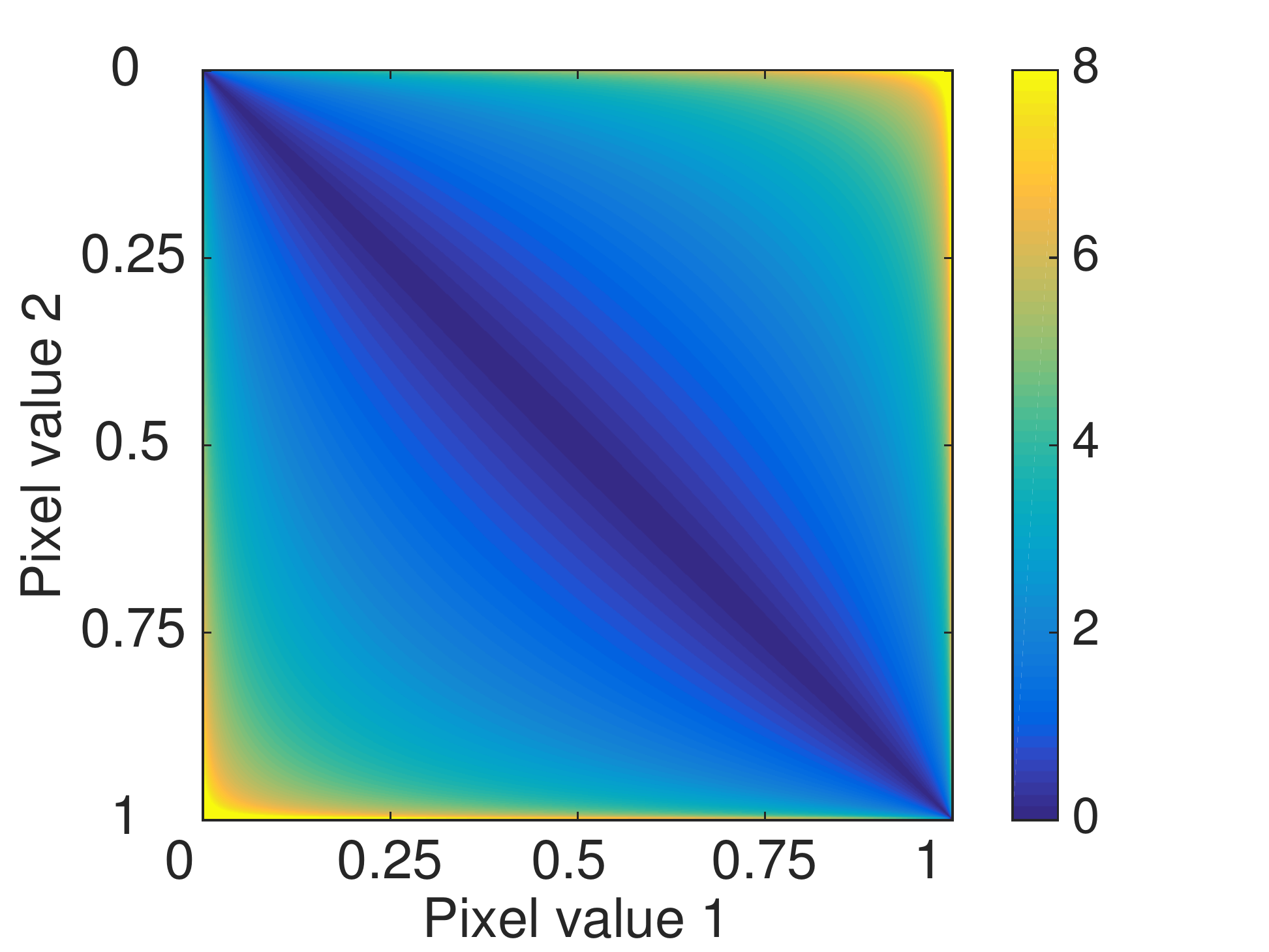}} ~
\subfloat[Proposed reparam. ($k=0.4$)]{\includegraphics[width=\figwidth]{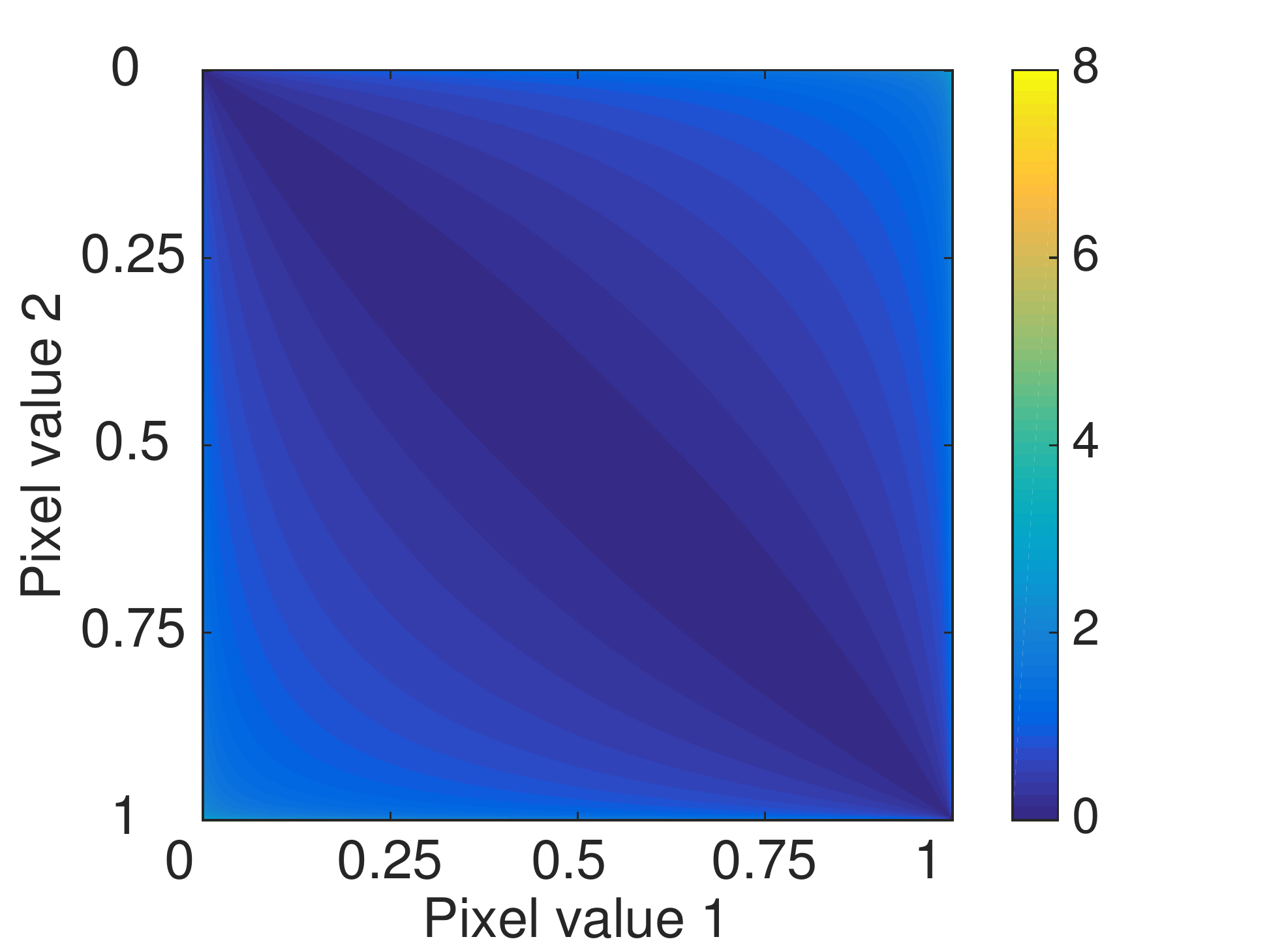}}
\captionsetup{singlelinecheck=off}
\caption[blah]{\textbf{Bernoulli implicit regularizers:} In these illustrations, each axis represents varying pixel values and the brightness shows the absolute difference $\left| [f^{-1}(x)]_1 - [f^{-1}(x)]_2 \right|$ for:

\vspace{1mm}
\begin{itemize}
\item[(a)] $f^{-1}(x) = x$ (using the $f(\theta) = \theta$ pixel intensity reparameterization)
\item[(b)] $f^{-1}(x) = \log(\frac{x}{1-x})$ (using the $f(\theta) = g^{-1}(\theta) = \frac{\exp \theta}{1+\exp \theta}$ natural reparameterization)
\item[(c)] $f^{-1}(x)$ as defined in \eqref{eq:bern_implicit} for $k=0.4$ (using the proposed reparameterization \eqref{eq:bern_f}).
\end{itemize}

\vspace{1mm}
Because the TV norm essentially takes image gradient differences, these illustrations give an intuition behind the implicit regularizers introduced by each reparameterization when using $\rho(\cdot) = \| \cdot \|_\mathrm{TV}$. We see that when pixel values are within $[0,k]$ and $(1-k,1]$ our proposed regularizer in (c) behaves like (b), but for pixel values in $(k,1-k]$ our proposed regularizer (c) behaves like (a); this is consistent with our piecewise construction of our proposed reparameterization.}
\label{fig:reg_bern}
\end{center}
\end{figure}

\subsection{Poisson noise}
\label{sec:poisson}

Poisson noise is a problem that is of wide interest in the image processing community, where it is also sometimes referred to as shot noise. This noise is prevalent in low-light imaging where the dominant noise source comes from photon counting; this arises in various fields of study where images are acquired in dark backgrounds (i.e. night vision) or are acquired with a short exposure time (i.e. biomedical imaging). In these photon-limited realms, the noise is modeled as being drawn from a Poisson distribution \cite{SPIRAL, snyder}, whose rate parameter represents the true underlying signal intensity, namely $y \sim \mathrm{Poisson}(x^*)$.

In the context of our generalized model from \eqref{eq:general_datafit}, applying Table \ref{tab:exp_fam} leaves us with
\begin{align}
\phi(f(\theta),y) = \sum_{i=1}^p \exp[g(f(\theta))]_i - y_i [g(f(\theta))]_i.
\label{eq:datafit_poisson}
\end{align}
One candidate $g(f(\theta))$ function comes from the generalization \eqref{eq:general_reparam}. For $\mathcal{F}$ as defined in \eqref{eq:family},
\begin{equation}
\begin{aligned}
[g(f(\theta))]_i \deq \begin{cases} k\theta_i + \log(\frac{1}{k}) - 1, & \theta_i \le \frac{1}{k} \\
\log(\theta_i), & \theta_i > \frac{1}{k} \end{cases},
\label{eq:reparam_poisson}
\end{aligned}
\end{equation}
where the image intensities are mapped back for some $k \in (0,\infty)$ via
\begin{equation}
\begin{aligned}
[f(\theta)]_i \deq \begin{cases} \exp(k\theta_i + \log(\frac{1}{k}) - 1), & \theta_i \le \frac{1}{k} \\
\theta_i, & \theta_i > \frac{1}{k} \end{cases},
\end{aligned}
\end{equation}
which is similar to the function proposed in \cite{hybridTV}.
In applying \eqref{eq:family} to the Poisson case, there are two degrees of freedom in $a$ and $b$. However, there is only one degree of freedom in the above expressions as a result of the pair $(a,b)$ satisfying continuity constraints. The family for reparameterizations can be formally rewritten as
\begin{equation}
\begin{aligned}
\mathcal{F} = \left\{ \vphantom{\begin{cases} a \\ a \end{cases}} \right.
f: \Theta \rightarrow \reals^p: \exists k>0 \ \text{s.t.} \ [f(\theta)]_i \deq \begin{cases} \exp(k\theta_i + \log(\frac{1}{k}) - 1), & \theta_i \le \frac{1}{k} \\
\theta_i, & \theta_i > \frac{1}{k} \end{cases}, \ i = 1,\cdots,p
\left. \vphantom{\begin{cases} a \\ a \end{cases}} \right \}.
\end{aligned}
\end{equation}

\subsubsection{Alternative Poisson reparameterization}

Another candidate reparameterization to explore is based on a function already proposed in the literature in \cite{SeegarBouchard2012}, where
\begin{align}
g(f(\theta)) \deq \log(k\log(1+\exp(\theta))),
\end{align}
where the image intensities are remapped via
\begin{align}
f(\theta) \deq k\log(1+\exp(\theta)),
\label{eq:logexp}
\end{align}
in what we call the ``$\log$-$\exp$'' parameterization. Though this function was motivated by a Poisson matrix completion problem, it was designed in order to address the unbounded curvature of the data-fit when estimating the rate parameter directly. We introduce the $k$ parameter in this work to describe a family of reparameterizations, which can be formally stated as
\begin{equation}
\begin{aligned}
\mathcal{F} = \Big\{ f: \Theta \rightarrow \reals^p: \exists k>0 \ \text{s.t.} \ f(\theta) \deq k\log(1+\exp(\theta)) \Big\}.
\end{aligned}
\end{equation}

\begin{framed}
Summarizing, in the Poisson case we propose the data-fit expression
\begin{align}
\phi(f(\theta),y) = \sum_{i=1}^p k\log(1+\exp(\theta_i)) - y_i \log(k\log(1+\exp(\theta_i)))
\label{eq:poisson_final}
\end{align}
for some $k \in (0,\infty)$.
\end{framed}

It is straightforward to show that the reparameterized data-fit term $\phi(f(\theta),y)$ we propose is convex and has Lipschitz gradients through algebra (convexity by showing that the diagonal Hessian matrix has diagonal elements $\frac{\partial^2 \phi(f(\theta),y)}{\partial \theta_i^2}$ that are nonnegative, and Lipschitz gradients by showing $\frac{\partial^2 \phi(f(\theta),y)}{\partial \theta_i^2}$'s are finite).

Our proposed $f(\theta)$ is an invertible function. Invertibility is not required in our setting, but it is useful in understanding the implicit regularization function $\rho(f^{-1}(x))$ seen in \eqref{eq:implicit_reg}. The inverse function is
\begin{align}
f^{-1}(x) = \log(\exp(\frac{x}{k})-1).
\label{eq:poiss_implicit}
\end{align}
For an intuition behind the regularization, Figure \ref{fig:reg_poiss} illustrates how the implicit regularizer $\| f^{-1}(x) \|_\mathrm{TV}$ behaves.

\setlength{\figwidth}{0.33\textwidth}
\begin{figure}[t]
\begin{center}
\subfloat[Intensity reparam.]{\includegraphics[width=\figwidth]{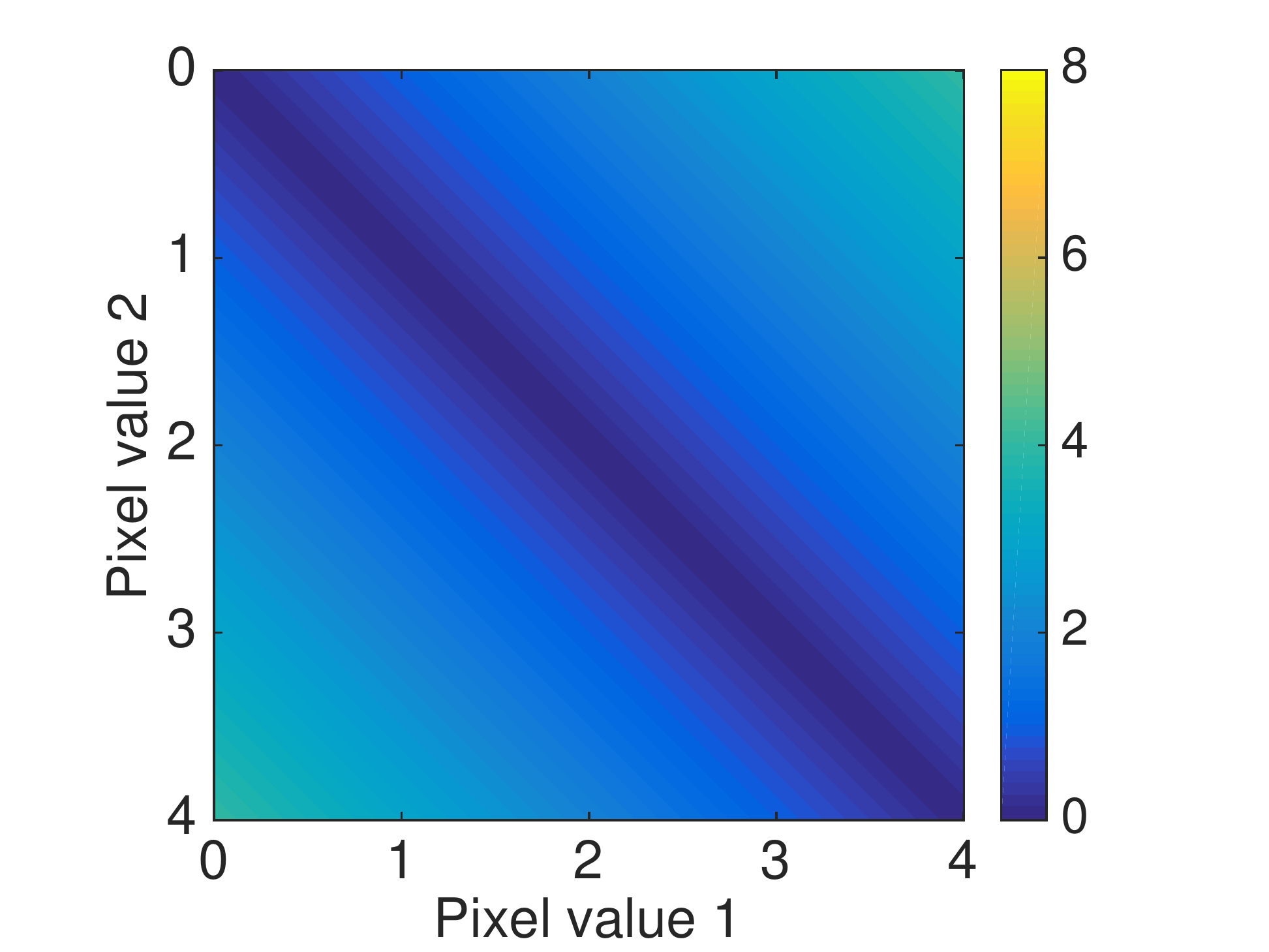}} ~
\subfloat[Natural reparam.]{\includegraphics[width=\figwidth]{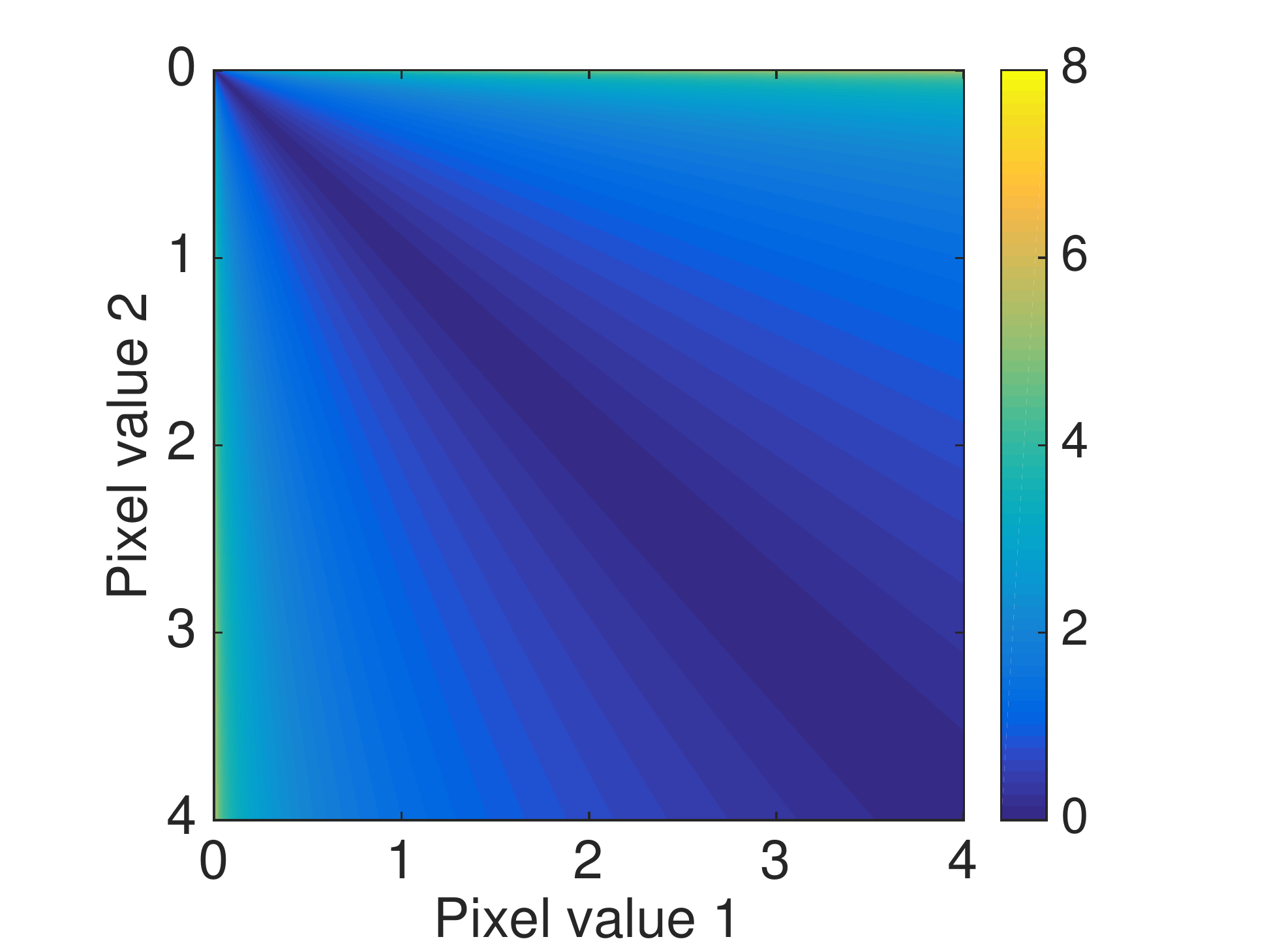}} ~
\subfloat[Proposed $\log$-$\exp$ ($k=2$)]{\includegraphics[width=\figwidth]{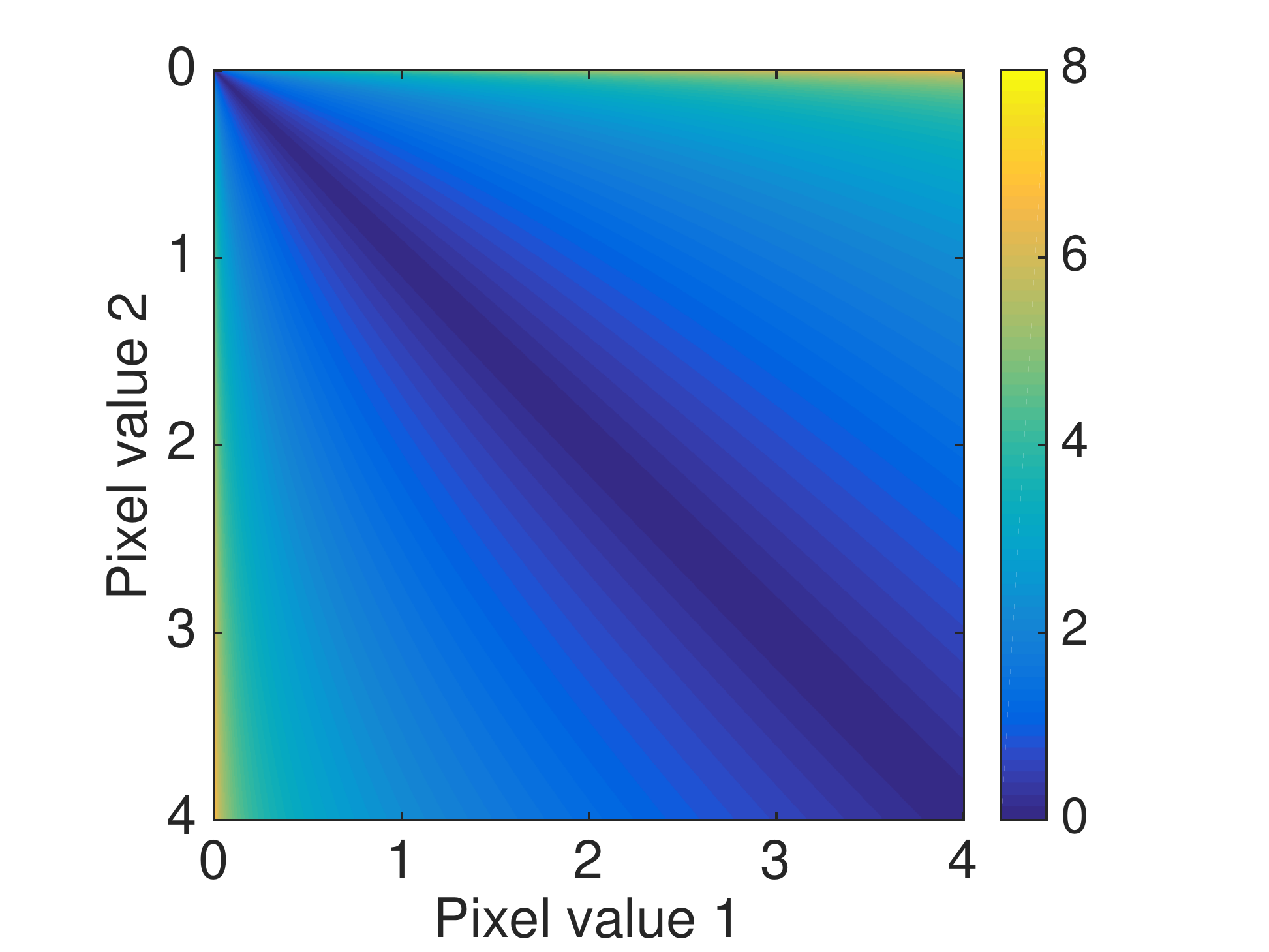}}
\captionsetup{singlelinecheck=off}
\caption[blah]{\textbf{Poisson implicit regularizers:} In these illustrations, each axis represents varying pixel values and the brightness shows the absolute difference $\left| [f^{-1}(x)]_1 - [f^{-1}(x)]_2 \right|$ for:

\vspace{1mm}
\begin{itemize}
\item[(a)] $f^{-1}(x) = x$ (using the $f(\theta) = \theta$ pixel intensity reparameterization)
\item[(b)] $f^{-1}(x) = \log(x)$ (using the $f(\theta) = g^{-1}(\theta) = \exp(\theta)$ natural reparameterization)
\item[(c)] $f^{-1}(x)$ as defined in \eqref{eq:poiss_implicit} for $k=2$ (using the proposed reparameterization \eqref{eq:logexp}).
\end{itemize}

\vspace{1mm}
Because the TV norm essentially takes image gradient differences, these illustrations give an intuition behind the implicit regularizers introduced by each reparameterization when using $\rho(\cdot) = \| \cdot \|_\mathrm{TV}$. We see that when pixel values are small the proposed regularizer in (c) behaves like (b), but for pixel values that are large our proposed regularizer (c) behaves like (a).}
\label{fig:reg_poiss}
\end{center}
\end{figure}

\section{Reparameterization for speckle noise}
\label{sec:gamma}

Speckle noise is found in many imaging applications, particularly in Synthetic Aperture Radar (SAR) imaging \cite{SARbook} and ultrasound \cite{ultrasoundSpeckle}. Speckle differs from other noise types in that it is multiplicative, where the observations obtained from an imaging system are the products of the true underlying intensity and draws from a Gamma distribution. The observation model can be written element-wise as follows:
\begin{align}
y_i = x_i n_i,
\end{align}
where $n$ is multiplicative noise that is known to be from a Gamma distribution with density
\begin{align}
p(n_i) = \frac{\beta}{\Gamma(\alpha)} n_i^{\alpha-1} \exp(-\beta n_i),
\end{align}
where $\alpha$ and $\beta$ are shape and rate parameters, respectively. In SAR imaging, for example, it is commonly assumed that $\alpha = \beta = S$, where $S$ is the number of ``looks" (or measurements) \cite{teuber2010}. Using this, we find that the negative log-likelihood $p(y|x)$ is then proportional to
\begin{align}
\phi(x,y) = \sum_{i=1}^p S \log x_i + S \frac{y_i}{x_i}.
\end{align}
Reparameterization results in the data-fit term
\begin{align}
\phi(f(\theta),y) = \sum_{i=1}^p S \log [f(\theta)]_i + S \frac{y_i}{[f(\theta)]_i},
\label{eq:datafit_gamma}
\end{align}
We first consider the reparameterization
\begin{align}
[f_1(\theta)]_i \deq 1/\theta_i
\label{eq:f1}
\end{align}
that yields the data-fit $\sum_{i=1}^p -\log \theta_i + y_i \theta_i$, which is convex; this function, however, has unbounded curvature as $\theta \rightarrow 0$. To remedy this, we can use an exponential transform via
\begin{align}
[f_2(\theta)]_i \deq 1/\exp(\theta_i)
\label{eq:f2}
\end{align}
to obtain $\sum_{i=1}^p -\theta_i + y_i\exp \theta_i$; this is also convex, but has unbounded curvature as $\theta \rightarrow \infty$. To address the unboundness of both possible solutions, we combine them through
\begin{equation}
\begin{aligned}
[f(\theta)]_i \deq \begin{cases} \frac{1}{k \exp(\theta_i)}, & \theta_i \le 0\\
\frac{1}{k (1+\theta_i)}, & \theta_i > 0 \end{cases},
\label{eq:reparam_gamma}
\end{aligned}
\end{equation}
where $k > 0$ affects the curvature, and allows us to define a family of reparameterizations over which to search.
Formally, the family this describes is
\begin{equation}
\begin{aligned}
\mathcal{F} = \left\{ f: \Theta \rightarrow \reals^p: \exists k>0 \ \text{s.t.} \ [f(\theta)]_i \deq \begin{cases} \frac{1}{k \exp(\theta_i)}, & \theta_i \le 0\\
\frac{1}{k (1+\theta_i)}, & \theta_i > 0 \end{cases}, \ i = 1,\cdots,p \right \}.
\end{aligned}
\end{equation}

\begin{framed}
Summarizing, in the speckle noise case we propose
\begin{align}
\phi(f(\theta),y) = \sum_{i=1}^p S \phi_i(f(\theta),y),
\label{eq:gamma_final}
\end{align}
where
\begin{align}
\phi_i(f(\theta),y) \deq \begin{cases} y_i k \exp(\theta_i) - \theta_i - \log k, & \theta_i \le 0\\
y_i k (1+\theta_i) - \log(1+\theta_i) - \log k, & \theta_i > 0 \end{cases}
\end{align}
for some $k \in (0,\infty)$.
\end{framed}

It is straightforward to show that the reparameterized data-fit term $\phi(f(\theta),y)$ we propose is convex and has Lipschitz gradients through algebra (convexity by showing that the diagonal Hessian matrix has diagonal elements $\frac{\partial^2 \phi(f(\theta),y)}{\partial \theta_i^2}$ that are nonnegative, and Lipschitz gradients by showing $\frac{\partial^2 \phi(f(\theta),y)}{\partial \theta_i^2}$'s are finite).

Our proposed $f(\theta)$ is an invertible function. Invertibility is not required in our setting, but it is useful in understanding the implicit regularization function $\rho(f^{-1}(x))$ seen in \eqref{eq:implicit_reg}. The inverse function is
\begin{align}
[f^{-1}(x)]_i = \begin{cases} \log(\frac{1}{kx_i}), & x_i \le \frac{1}{k}\\
\frac{1}{kx_i}-1, &x_i > \frac{1}{k} \end{cases}.
\label{eq:gamma_implicit}
\end{align}
For an intuition behind the regularization, Figure \ref{fig:reg_gamma} illustrates how the implicit regularizer $\| f^{-1}(x) \|_\mathrm{TV}$ behaves.

\setlength{\figwidth}{0.33\textwidth}
\begin{figure}[t]
\begin{center}
\subfloat[$f_1(\theta)$ reparam.]{\includegraphics[width=\figwidth]{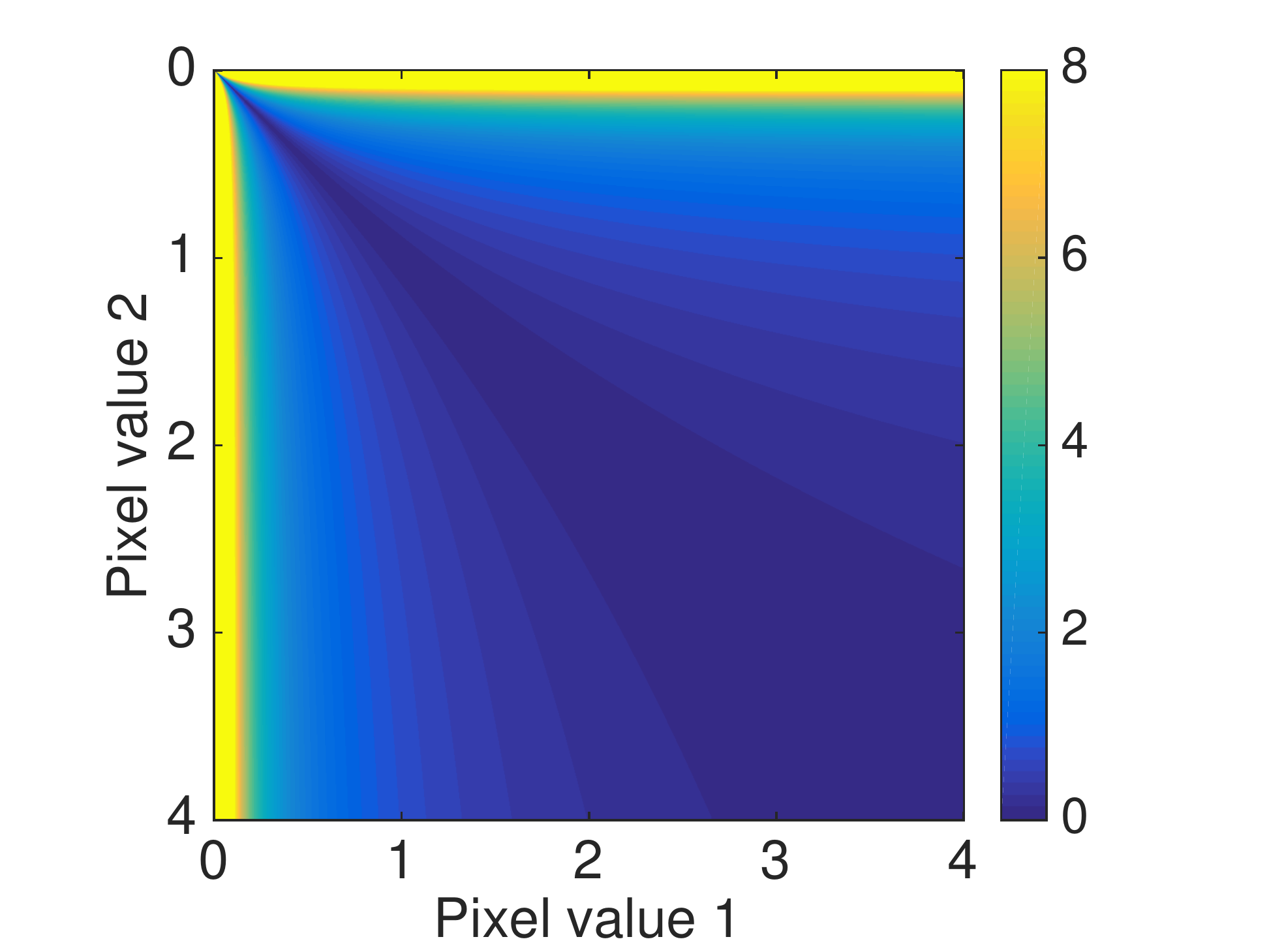}} ~
\subfloat[$f_2(\theta)$ reparam.]{\includegraphics[width=\figwidth]{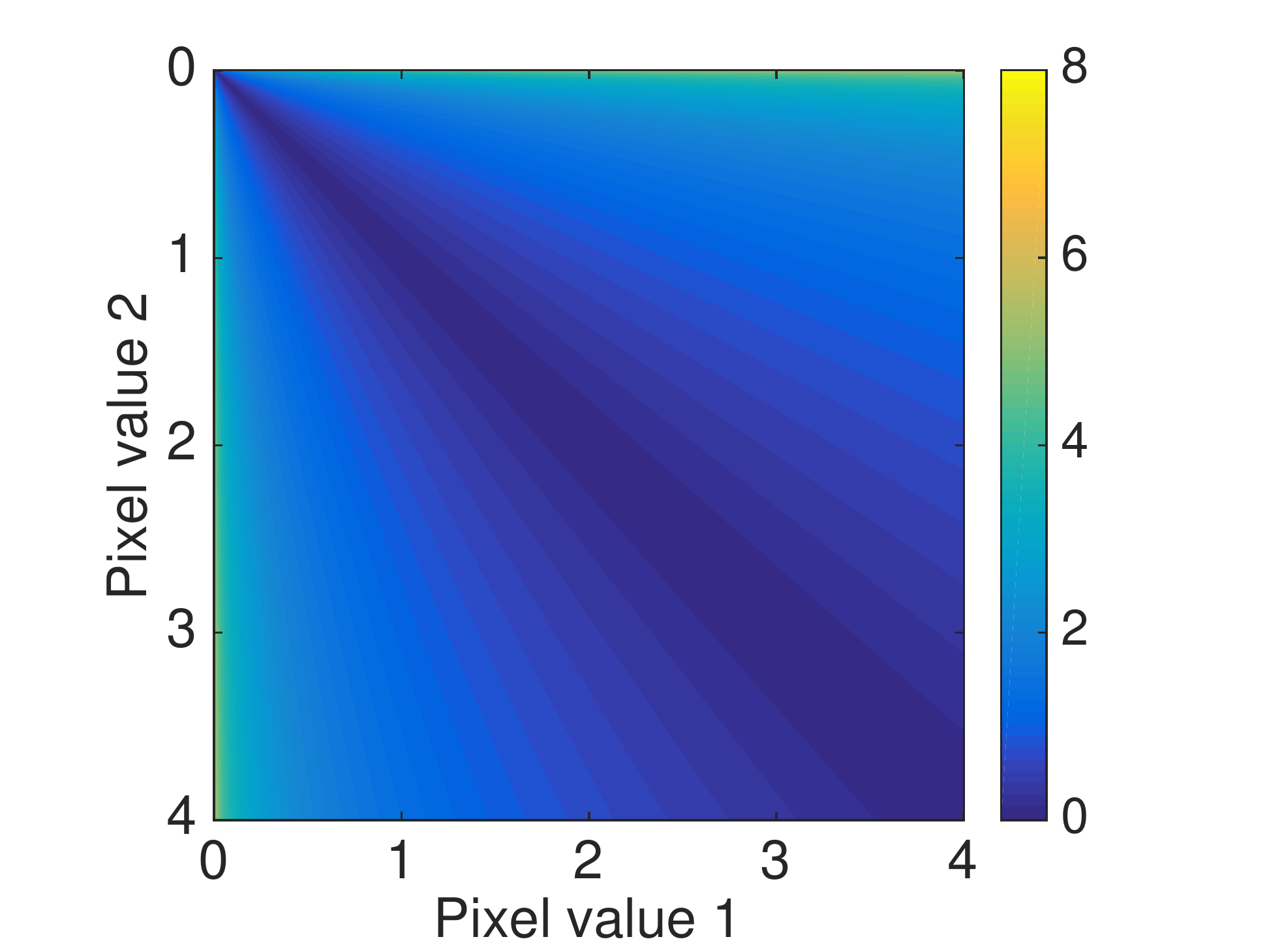}} ~
\subfloat[Proposed reparam. ($k=2$)]{\includegraphics[width=\figwidth]{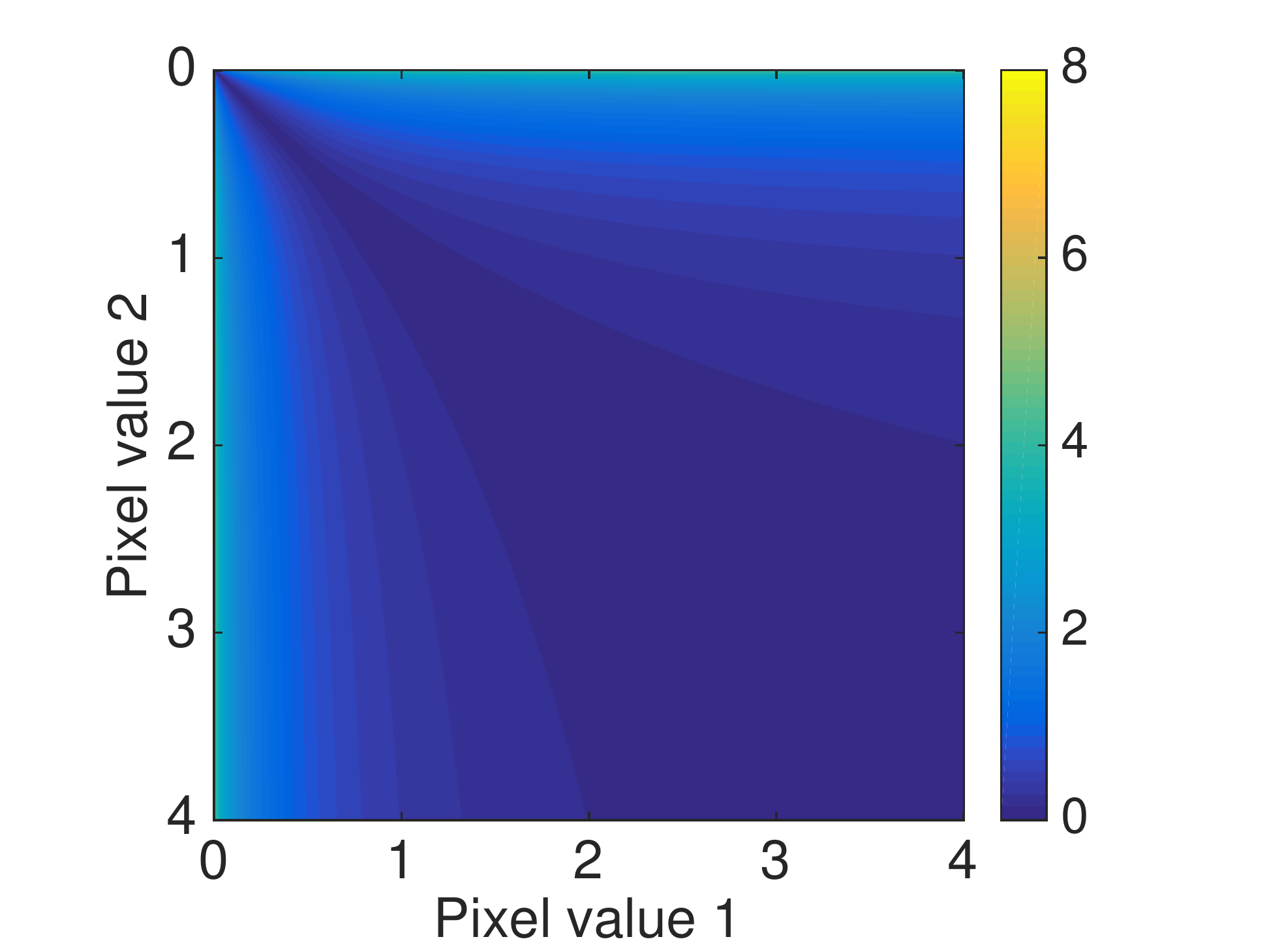}}
\captionsetup{singlelinecheck=off}
\caption[blah]{\textbf{Speckle noise implicit regularizers:} In these illustrations, each axis represents varying pixel values and the brightness shows the absolute difference $\left| [f^{-1}(x)]_1 - [f^{-1}(x)]_2 \right|$ for:

\vspace{1mm}
\begin{itemize}
\item[(a)] $[f^{-1}(x)]_i = \frac{1}{x_i}$ (using the $[f_1(\theta)]_i = \frac{1}{\theta_i}$ reparameterization)
\item[(b)] $f^{-1}(x) = -\log(x)$ (using the $[f_2(\theta)]_i = \frac{1}{\exp(\theta_i)}$ reparameterization)
\item[(c)] $f^{-1}(x)$ as defined in \eqref{eq:gamma_implicit} for $k=2$ (using the proposed reparameterization \eqref{eq:reparam_gamma}).
\end{itemize}

\vspace{1mm}
Because the TV norm essentially takes image gradient differences, these illustrations give an intuition behind the implicit regularizers introduced by each reparameterization when using $\rho(\cdot) = \| \cdot \|_\mathrm{TV}$. We see that for small pixel values within $[0,\frac{1}{k}]$ the proposed regularizer in (c) behaves like (b), but for large pixel values within $(\frac{1}{k},+\infty)$ our proposed regularizer (c) behaves like (a); this is consistent with the piecewise construction of our proposed reparameterization.}
\label{fig:reg_gamma}
\end{center}
\end{figure}

\section{Numerical experiments}
\label{sec:experiments}

As described in Section \ref{sec:approach}, we seek the reparameterization $f$ whose resulting estimate $\hat{x}$ minimizes the average normalized RMSE ($\| \hat{x} - x^* \|_2 / \|x^*\|_2$) over a family $\mathcal{F}$
for a bank of images for a given noise level. In other words, we are interested in finding the reparameterization $f$ that minimizes the empirical risk $\hat{R}(f,\rho)$ as in \eqref{eq:relative_risk}, or
\begin{align}
\hat{f} = \underset{f \in \mathcal{F}}{\argmin} \ \hat{R}(f,\rho).
\label{eq:risk_numexp}
\end{align}

To be clear about the experimental results to follow (and to tersely sum up our discussion in previous sections), we obtain the estimates $\hat{x}_{f,\rho,\tau}^{(i,j)}$ for the $i^{th}$ image and $j^{th}$ noise realization by using: data-fit \eqref{eq:bernoulli_final} for Bernoulli noise, $\log$-$\exp$ data-fit \eqref{eq:poisson_final} for Poisson noise, and data-fit \eqref{eq:gamma_final} for speckle noise. To showcase the value of our reparameterizations, we compare our experimental results to a baseline denoising that minimizes the negative log-likelihood with a regularizer applied directly on the image intensities, namely,
\begin{align}
\tilde{x}_{\tau} &= \argmin_x \ \{ \Phi(x,y) = \phi(x,y) + \tau \rho(x) \},
\label{eq:baseline}
\end{align}
which is parameterized by the tuning weight $\tau$.

In the examples to follow, we use the TV regularizer as defined in \eqref{eq:TV} for $\rho$. We note that all $\tau$ weighting values for the reparameterization and baseline estimates in the following experiments were chosen clairvoyantly, as to minimize each realization's RMSE.

For our proximal gradient method, we adopt SpaRSA \cite{sparsa} to compute the minimizer of $\phi(f(\theta),y^{(i,j)})+\tau \| \theta \|_{\rm TV}$, for the $i^{th}$ image and its $j^{th}$ noise realization. The main idea of the method is to iteratively perform \eqref{eq:prox}, where $s_t \triangleq \theta^t - \frac{1}{\alpha_t} \nabla \phi(f(\theta^t),y^{(i,j)})$ is a step in the direction of steepest descent and $\alpha_t$ is a step-size that changes at each iteration.
For a TV regularizer, \eqref{eq:prox} can be solved with the Fast Gradient Projection (FGP) algorithm \cite{FGP}. We defer more detailed descriptions of both SpaRSA and FGP to their original papers, as they are not the main focus of this work but are used for experimental validation.

\subsection{Shepp-Logan phantom}

As a first example to demonstrate the value of reparameterizing, we show example reconstructions with the Shepp-Logan phantom ($512 \times 512$) for each of the noise cases considered. We note that for Bernoulli noise ground truth was scaled to have intensities between 0 and 1, for Poisson noise ground truth was scaled to be between 0 and 5, and for speckle noise ground truth was scaled to be between 0 and 5 (with 3 looks in the noise generation). Figure \ref{fig:shepp} shows reconstruction results for the baseline reconstruction \eqref{eq:baseline} and the best reparameterized reconstruction, and also lists the RMSE averaged over 20 noise realizations for each estimate; we see a significant improvement in terms RMSE and visual quality when reparameterizing. Figure \ref{fig:shepp_zoom} shows the same reconstruction results on a zoomed-in view for more detail.

\setlength{\figwidth}{0.2\textwidth}
\begin{figure}[ht]
\begin{center}
\subfloat[Truth]{\includegraphics[width=\figwidth]{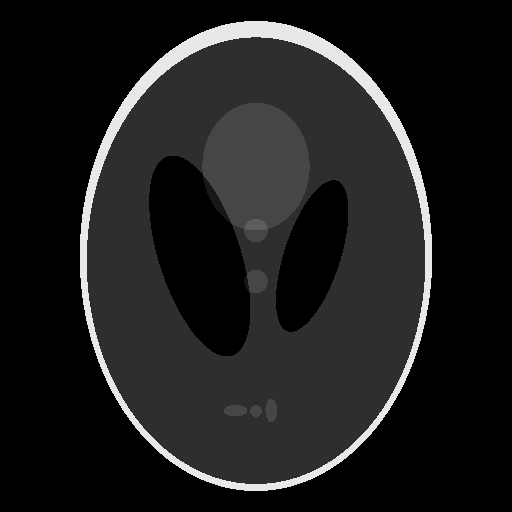}} ~
\subfloat[Noisy observations (Bernoulli)]{\includegraphics[width=\figwidth]{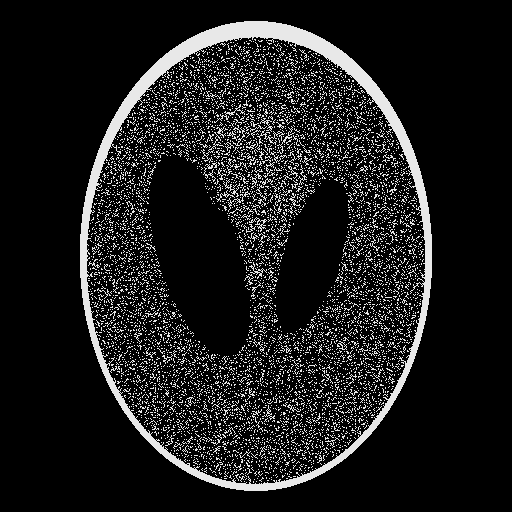}} ~
\subfloat[Baseline recovery \bf{(RMSE: 0.1554)}]{\includegraphics[width=\figwidth]{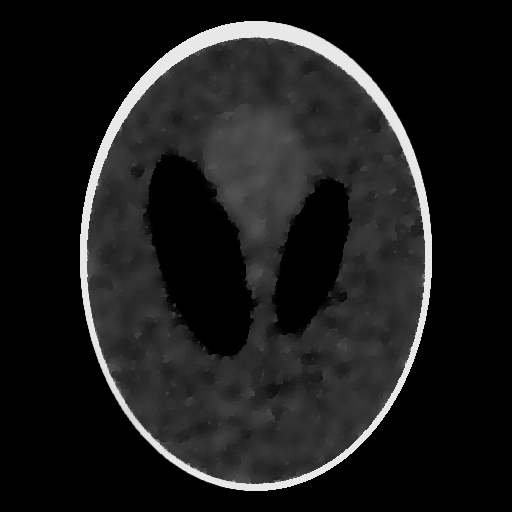}} ~
\subfloat[Reparameterized recovery, $k=0.05$ \bf{(RMSE: 0.1249)}]{\includegraphics[width=\figwidth]{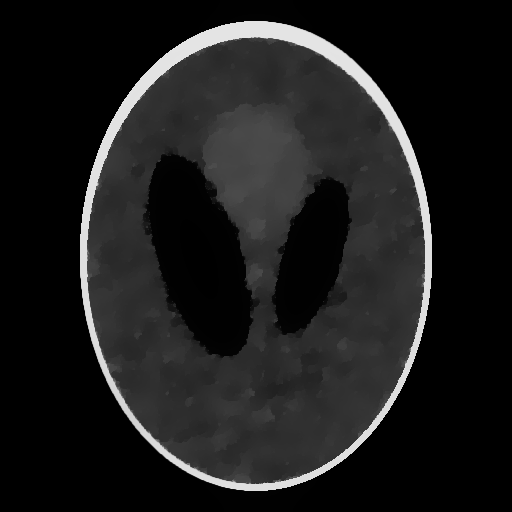}} \\
\subfloat[Truth]{\includegraphics[width=\figwidth]{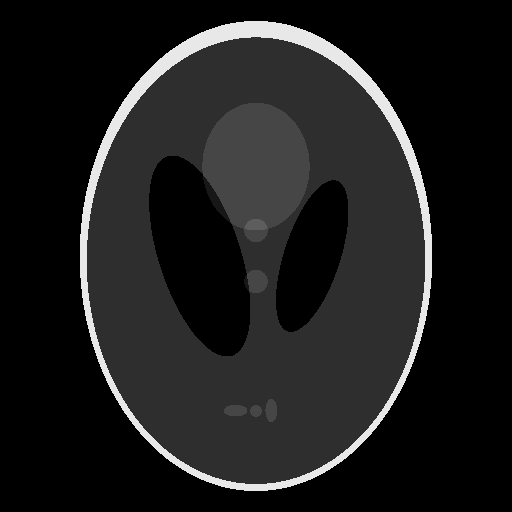}} ~
\subfloat[Noisy observations (Poisson)]{\includegraphics[width=\figwidth]{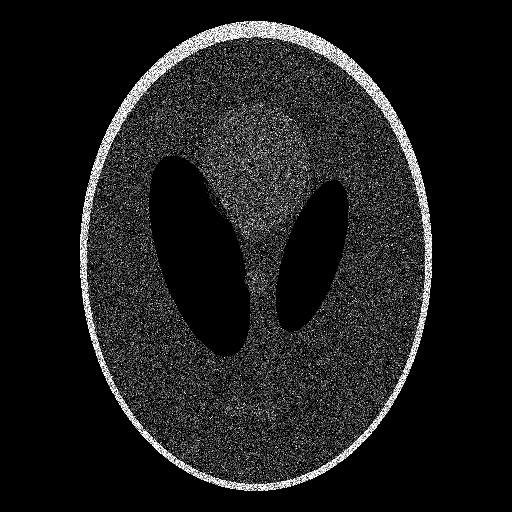}} ~
\subfloat[Baseline recovery \bf{(RMSE: 0.2074)}]{\includegraphics[width=\figwidth]{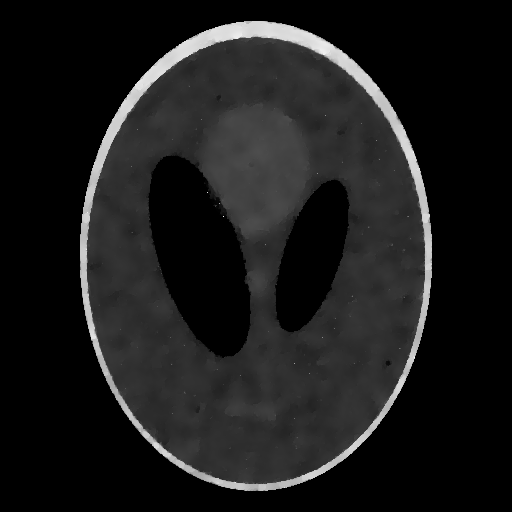}} ~
\subfloat[Reparameterized recovery, $k=4.0$ \bf{(RMSE: 0.1442)}]{\includegraphics[width=\figwidth]{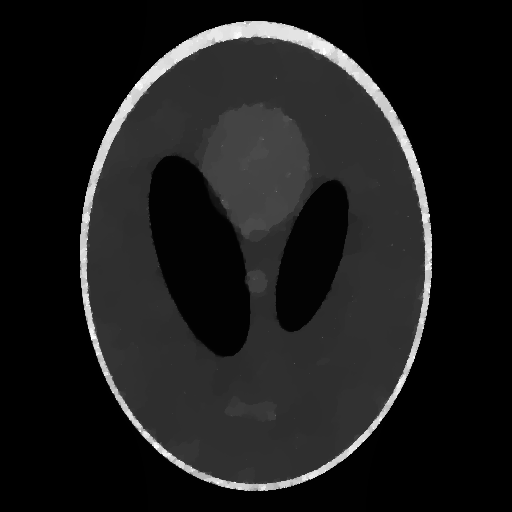}} \\
\subfloat[Truth]{\includegraphics[width=\figwidth]{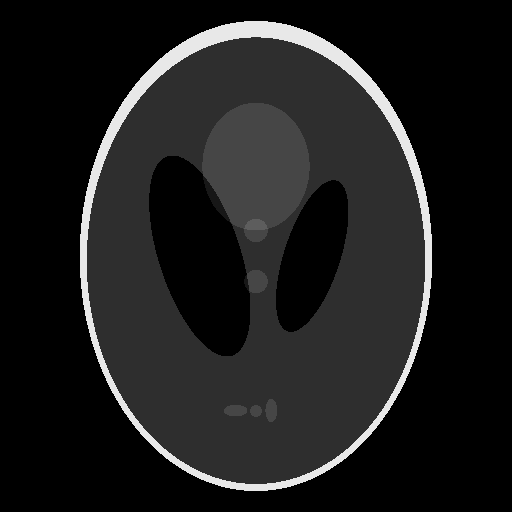}} ~
\subfloat[Noisy observations (Speckle)]{\includegraphics[width=\figwidth]{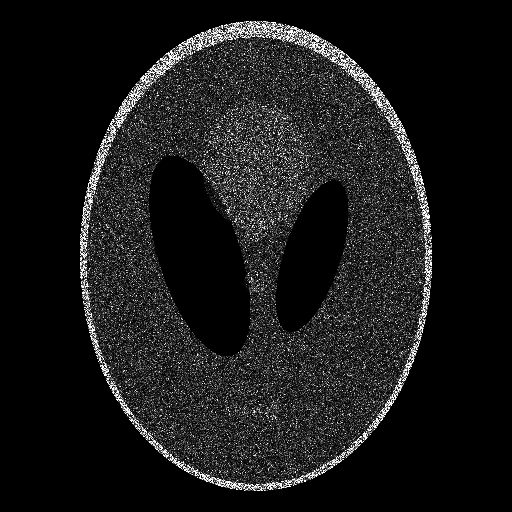}} ~
\subfloat[Baseline recovery \bf{(RMSE: 0.6483)}]{\includegraphics[width=\figwidth]{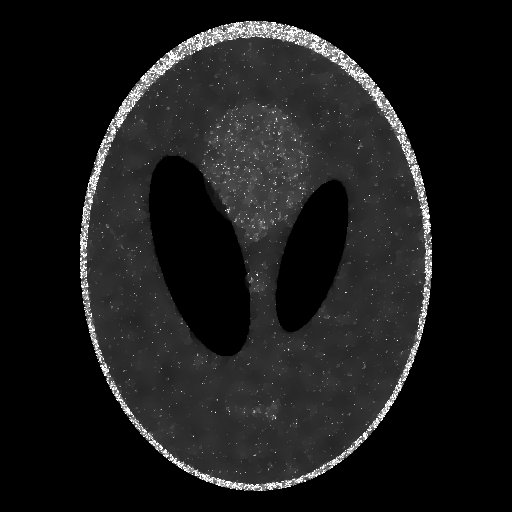}} ~
\subfloat[Reparameterized recovery, $k=2.0$ \bf{(RMSE: 0.3229)}]{\includegraphics[width=\figwidth]{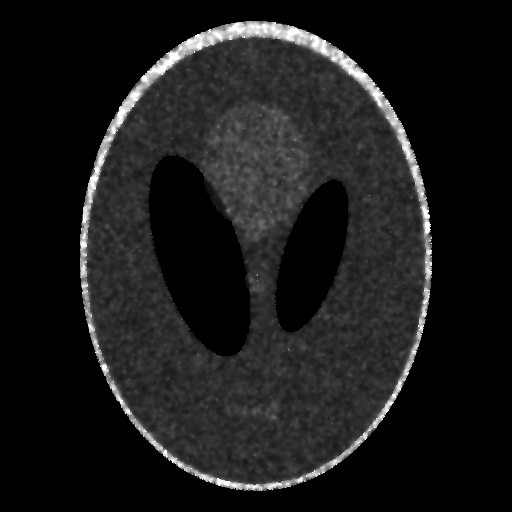}}~
\caption{\small \textbf{Shepp-Logan phantom experiments:} Denoising results using Bernoulli noise (first row), Poisson noise (second row), and speckle noise (third row). We show the ground truth, noisy observations, the baseline recovery \protect\eqref{eq:baseline}, and the best reparameterized recovery. The RMSEs of our reconstructions (averaged over 20 noise realizations) are also listed.}
\label{fig:shepp}
\end{center}
\end{figure}

\setlength{\figwidth}{0.2\textwidth}
\begin{figure}[ht]
\begin{center}
\subfloat[Truth]{\includegraphics[width=\figwidth]{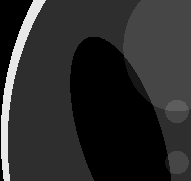}} ~
\subfloat[Noisy observations (Bernoulli)]{\includegraphics[width=\figwidth]{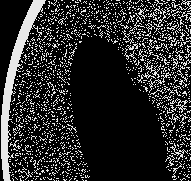}} ~
\subfloat[Baseline recovery \bf{(RMSE: 0.1554)}]{\includegraphics[width=\figwidth]{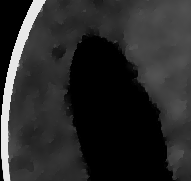}} ~
\subfloat[Reparameterized recovery, $k=0.05$ \bf{(RMSE: 0.1249)}]{\includegraphics[width=\figwidth]{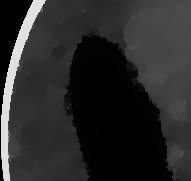}} \\
\subfloat[Truth]{\includegraphics[width=\figwidth]{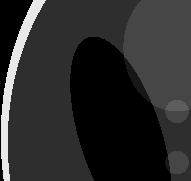}} ~
\subfloat[Noisy observations (Poisson)]{\includegraphics[width=\figwidth]{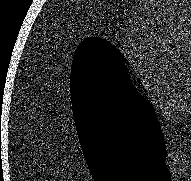}} ~
\subfloat[Baseline recovery \bf{(RMSE: 0.2074)}]{\includegraphics[width=\figwidth]{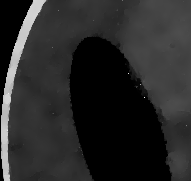}} ~
\subfloat[Reparameterized recovery, $k=4.0$ \bf{(RMSE: 0.1442)}]{\includegraphics[width=\figwidth]{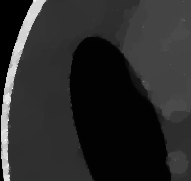}} \\
\subfloat[Truth]{\includegraphics[width=\figwidth]{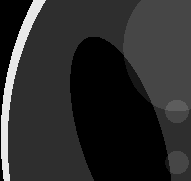}} ~
\subfloat[Noisy observations (Speckle)]{\includegraphics[width=\figwidth]{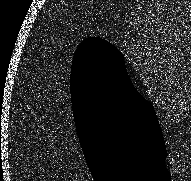}} ~
\subfloat[Baseline recovery \bf{(RMSE: 0.6483)}]{\includegraphics[width=\figwidth]{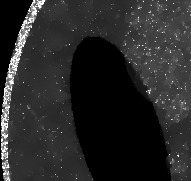}} ~
\subfloat[Reparameterized recovery, $k=2.0$ \bf{(RMSE: 0.3229)}]{\includegraphics[width=\figwidth]{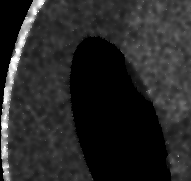}}~
\caption{\small \textbf{Shepp-Logan phantom experiments (zoomed-in view):} Denoising results using Bernoulli noise (first row), Poisson noise (second row), and speckle noise (third row). We show the ground truth, noisy observations, the baseline recovery \protect\eqref{eq:baseline}, and the best reparameterized recovery. The RMSEs of our reconstructions (averaged over 20 noise realizations) are also listed; note that the RMSEs are for the entire image and not for the zoomed-in parts.}
\label{fig:shepp_zoom}
\end{center}
\end{figure}

\subsection{Natural image experiments}
The Shepp-Logan experiments show value in reparameterizing in terms of RMSE and visual image quality. We now more extensively analyze reparameterizations through comparing reconstruction results on natural images to see which choice of reparameterization minimizes the empirical risk \eqref{eq:risk_numexp}. The bank of images considered consists of $n=10$ images from \cite{Dabov09bm3dimage}. These are commonly seen in the image processing literature, and include: Barbara (512$\times$512), Boat (512$\times$512), Cameraman (256$\times$256), Couple (512$\times$512), Fingerprint (512$\times$512), Hill (512$\times$512), House (256$\times$256), Lena (512$\times$512), Man (512$\times$512), and Peppers (256$\times$256). There are $m=20$ noisy realizations for each test image. For the Bernoulli case, the true images were scaled to have intensities on $[0,1]$, and for the Poisson case, the images were scaled to have nonnegative intensities of mean 5 and 10 per pixel. For the speckle noise experiments, the true images were scaled to have nonnegative intensities of mean 3 and 5 per pixel, with $S=3$ and $5$ looks.

\subsubsection{Exponential family experiments}

Figure \ref{fig:bern_poiss_plots} shows plots of how the relative risk $\hat{R}(f,\rho)$ with respect to $k$ for the Bernoulli and Poisson cases, using the proposed \eqref{eq:bernoulli_final} and $\log$-$\exp$ \eqref{eq:poisson_final}, respectively. These curved are averaged over the 10 images with 20 noise realizations each and indicate optimal settings of the parameter $k$ in each case. The data for these plots can be found in Tables \ref{table:table_bernoulli}, \ref{table:table_poisson_mean5}, and \ref{table:table_poisson_mean10} in the appendices, where the averaged RMSEs for each image's noise realizations are listed; the last row of the tables corresponds to the $\hat{R}$ plots in Figure \ref{fig:bern_poiss_plots}. We see that the baseline comparison performs worse than the reparameterizations considered. 

The appendices also show example reconstructions. Figure \ref{fig:bernoulli_images} shows example reconstructions from one noise realization of the Hill and House data sets under Bernoulli noise. Figure \ref{fig:poisson_images_mean5} shows example reconstructions from the Poisson realizations for the Barbara and Couple images with mean 5. Figure \ref{fig:poisson_images_mean10} shows example reconstructions from the Poisson realizations for the Lena and Man images with mean 10.

We note that for these noise examples, the baseline negative log-likelihood with respect to image intensities
is convex, which fits well into the SpaRSA framework.

\setlength{\figwidth}{0.33\textwidth}
\begin{figure}[t]
\begin{center}
\subfloat[Bernoulli]{\includegraphics[width=\figwidth]{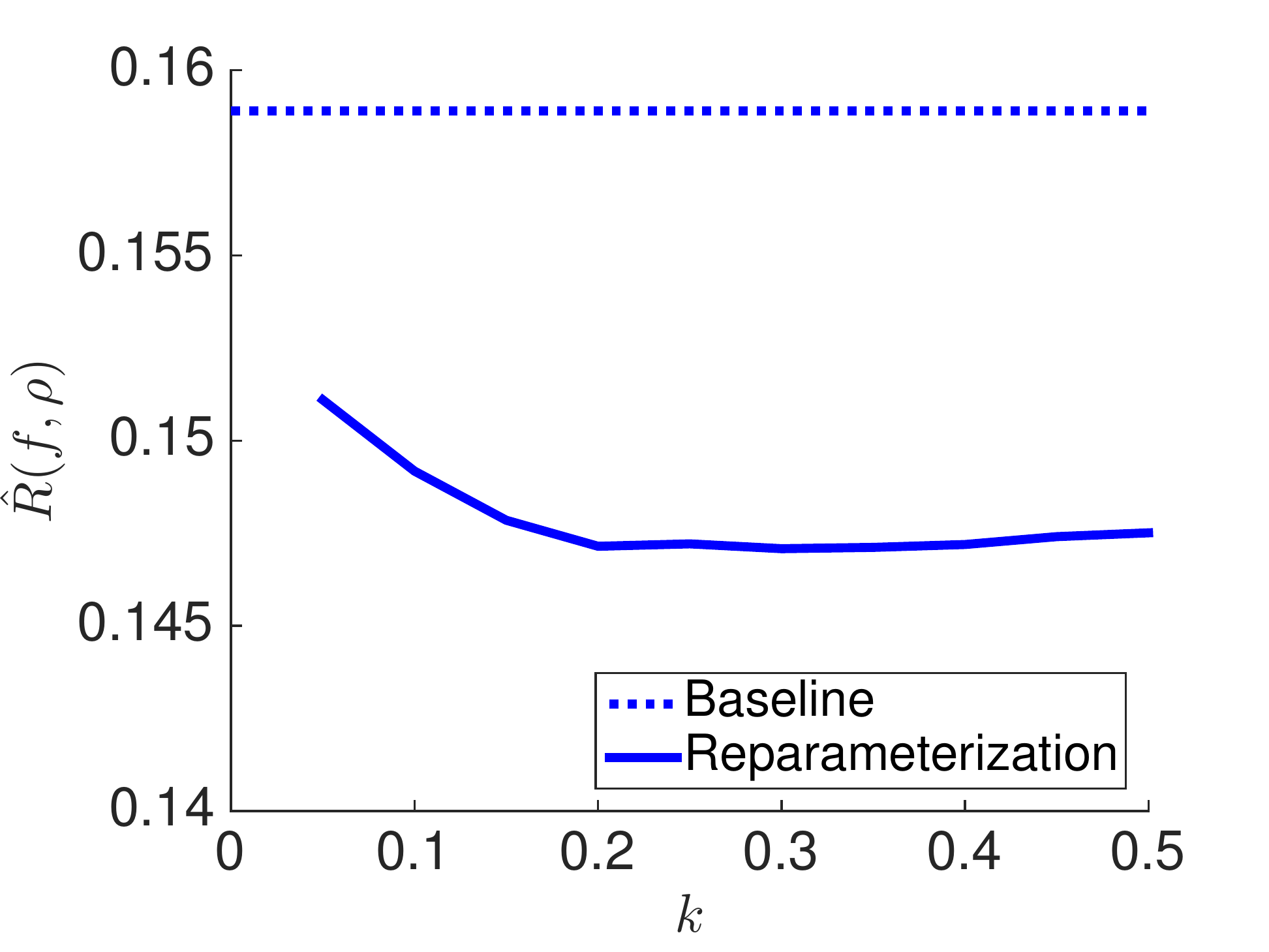} \label{fig:bern_plots}} ~
\subfloat[Poisson (mean 5)]{\includegraphics[width=\figwidth]{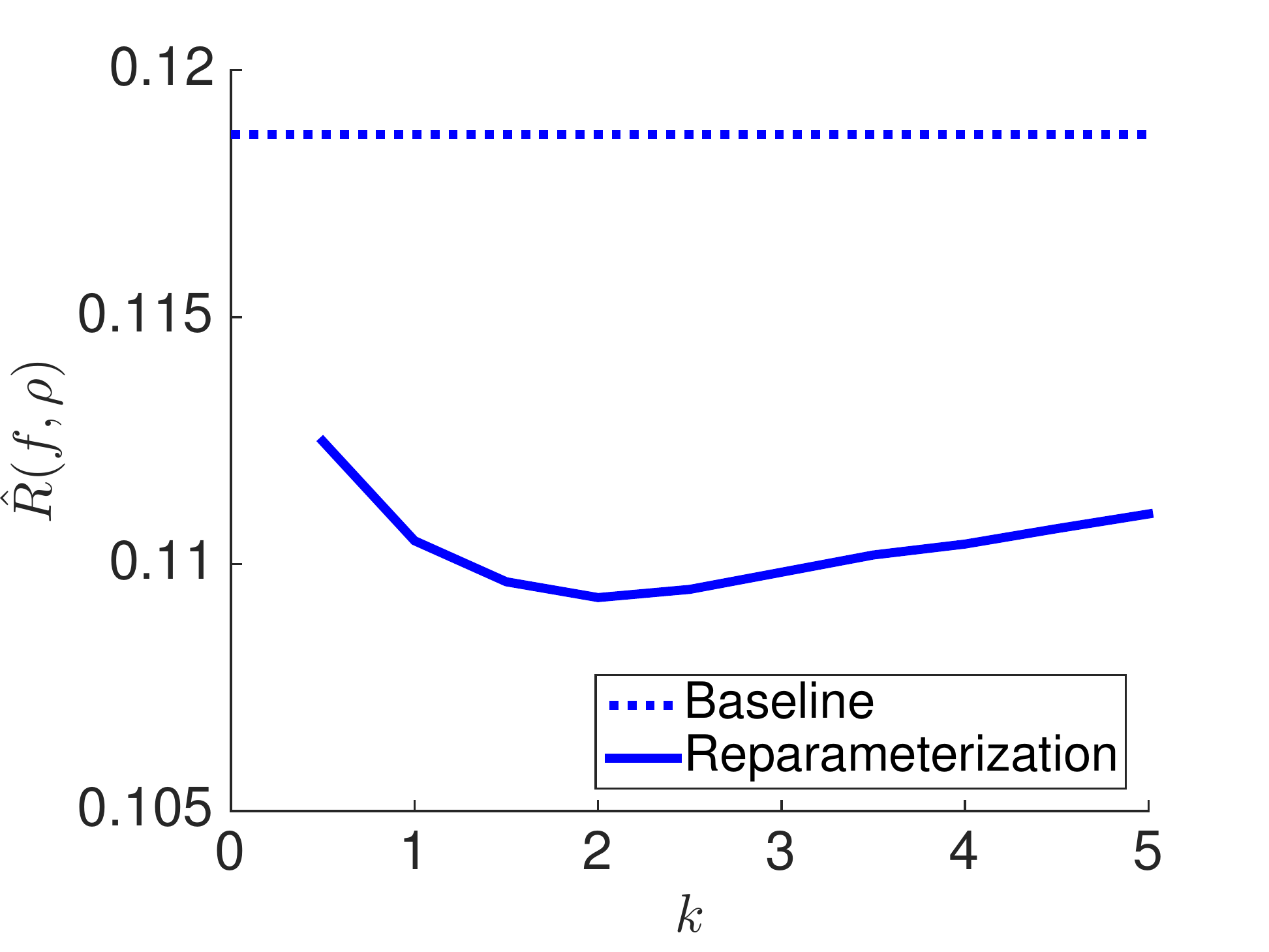} \label{fig:poiss_mean5_plots}} ~
\subfloat[Poisson (mean 10)]{\includegraphics[width=\figwidth]{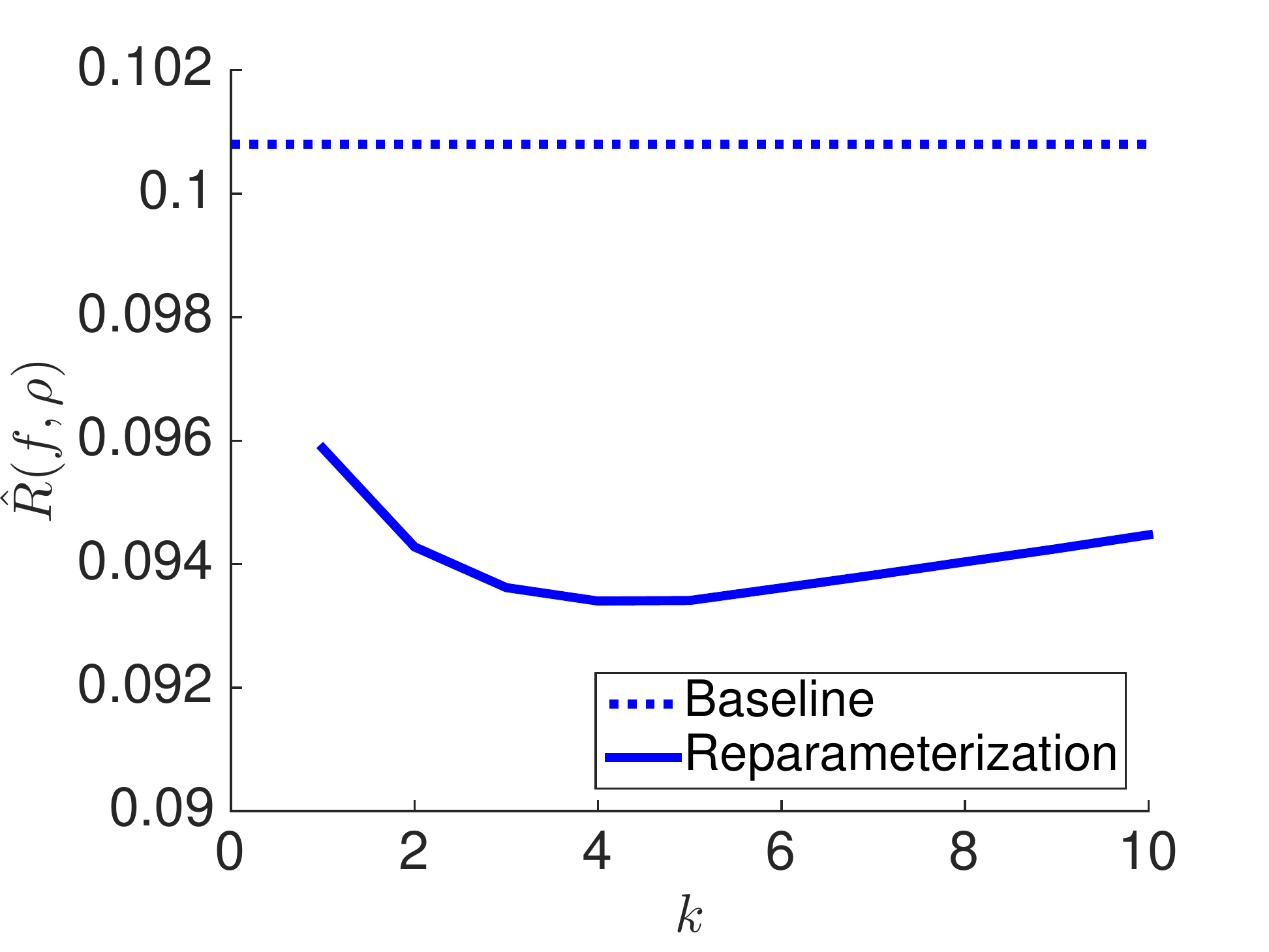} \label{fig:poiss_mean10_plots}} ~
\caption{\textbf{Bernoulli and Poisson noise $\hat{R}(f,\rho)$ plots:} These plots show how the relative risk $\hat{R}(f,\rho)$ for the proposed reparameterizations (solid lines) varies as the parameter $k$ changes for Bernoulli and Poisson noise, using the proposed \protect\eqref{eq:bernoulli_final} and $\log$-$\exp$ \protect\eqref{eq:poisson_final}, respectively. Also shown are $\hat{R}(f,\rho)$'s for the baseline reconstructions \protect\eqref{eq:baseline} (dashed lines) as comparisons. These plots are averages over 10 images with 20 noise realizations each; the data used to generate these plots can be found in Tables \ref{table:table_bernoulli}, \ref{table:table_poisson_mean5}, and \ref{table:table_poisson_mean10} in the appendices.}
\label{fig:bern_poiss_plots}
\end{center}
\end{figure}

\subsubsection{Speckle noise experiments}

Figure \ref{fig:gamma_plots} shows how $\hat{R}(f,\rho)$ changes as the reparameterization functions change with respect to $k$ for speckle noise using the proposed \eqref{eq:gamma_final}. Averaged over the 10 images with 20 noise realizations each, we see that as $k$ grows, $\hat{R}$ initially decreases but then stabilizes for large $k$. We remark that we search over small values of $k$ (between 0 and 2) for practical reasons, as large $k$ can cause numerical instabilities. Table \ref{table:table_gamma} in the appendices lists the RMSE averages of the baseline comparison and $k = 2$ reparameterization for each image considered; we only show the values at $k = 2$ because of the observed leveling off behavior.

Figure \ref{fig:gamma_images_mean3} shows example reconstructions from one noise realization of speckle noise of the Peppers ($S = 3$) and Cameraman ($S = 5$) data sets with mean 3. Figure \ref{fig:gamma_images_mean5} shows example reconstructions from one noise realization of speckle noise of the Fingerprint ($S=3$) and Boat ($S=5$) data sets with mean 5.

We note that the baseline negative log-likelihood data-fit for the speckle noise is not convex, but we can still apply SpaRSA to find a minimizer.

\setlength{\figwidth}{0.33\textwidth}
\begin{figure}[t]
\begin{center}
\subfloat[Mean 3]{\includegraphics[width=\figwidth]{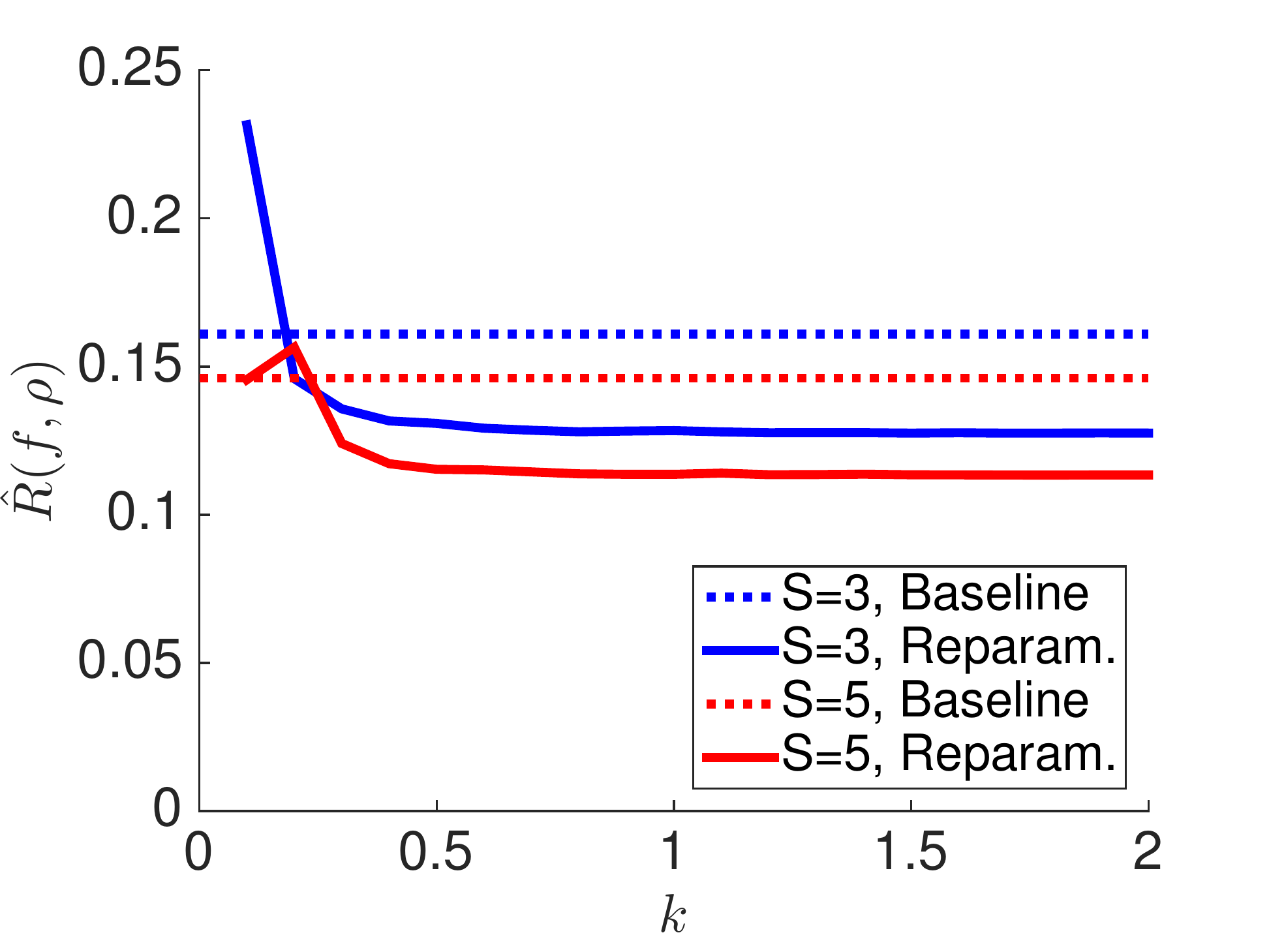}} ~
\subfloat[Mean 5]{\includegraphics[width=\figwidth]{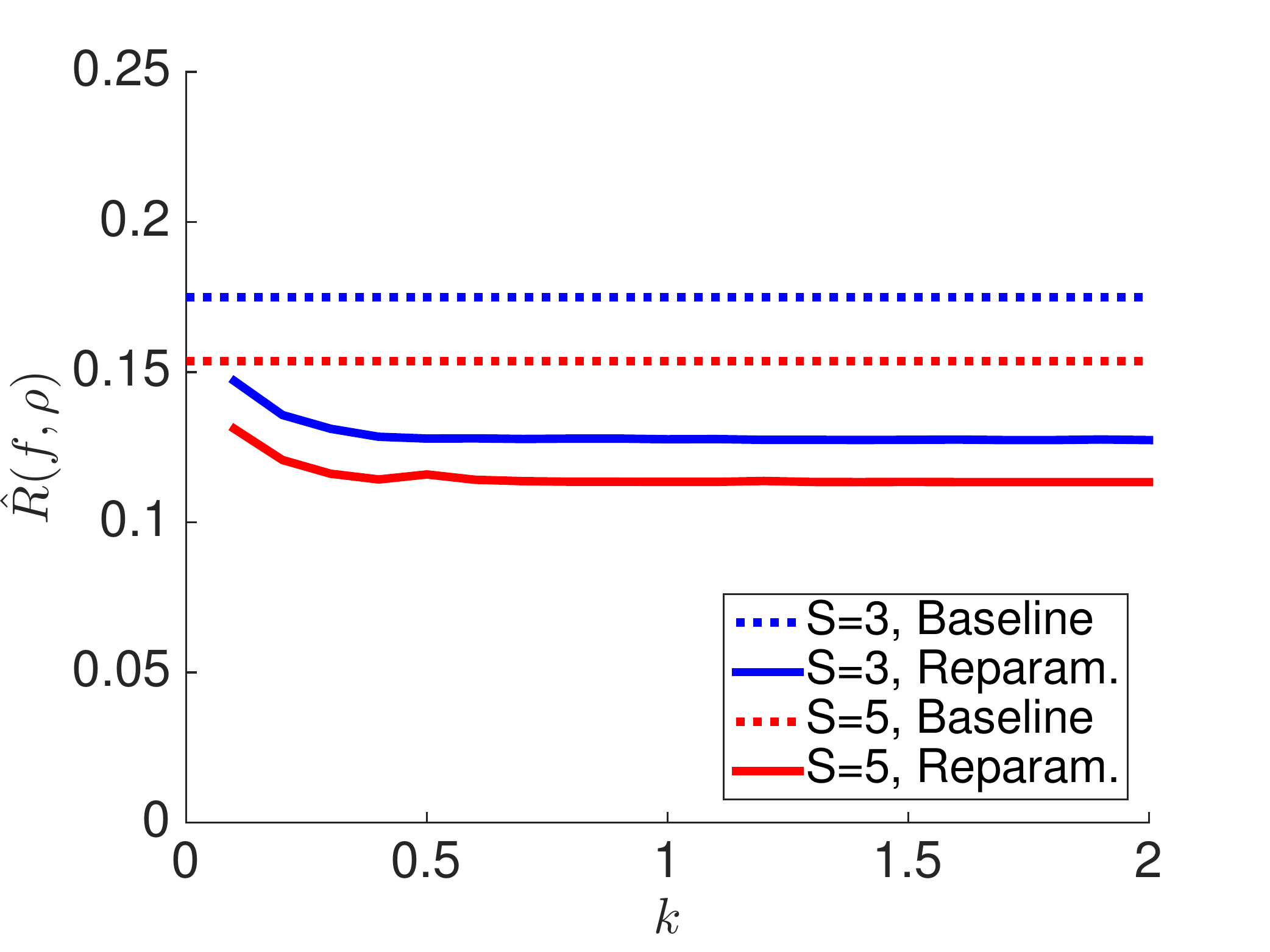}} ~
\caption{\textbf{Speckle noise $\hat{R}(f,\rho)$ plots:} These plots show how the relative risk $\hat{R}(f,\rho)$ for the proposed reparameterization (solid lines) varies as the parameter $k$ changes for various combinations of mean image intensity and number of looks $S$, using the proposed \protect\eqref{eq:gamma_final}. We see that $\hat{R}(f,\rho)$ levels off as $k$ increases. Also shown are $\hat{R}(f,\rho)$'s for the baseline reconstructions \protect\eqref{eq:baseline} (dashed lines) as comparisons. These plots are averages over 10 images with 20 noise realizations each. Table \ref{table:table_gamma} in the appendices shows the RMSE values for the baseline comparison and the $k=2$ reconstruction for each image considered.}
\label{fig:gamma_plots}
\end{center}
\end{figure}

\subsubsection{Computational-statistical tradeoff}

We briefly revisit the idea of a computational-statistical tradeoff discussed in Section \ref{sec:compstat}. Figure \ref{fig:compstat} illustrates one example from Poisson noise (specifically, we are looking at results from the Couple image with mean 10). The plots show how RMSE and the number of SpaRSA iterations behave with the parameter $k$. We see that RMSE vs. $k$ has a dip and that Iterations vs. $k$ has a generally increasing behavior. Plotting RMSE vs.\ iterations shows us how well we can expect to perform given a limited computational budget. 

\setlength{\figwidth}{0.3\textwidth}
\begin{figure}[t]
\begin{center}
\subfloat[RMSE vs. $k$]{\includegraphics[width=\figwidth]{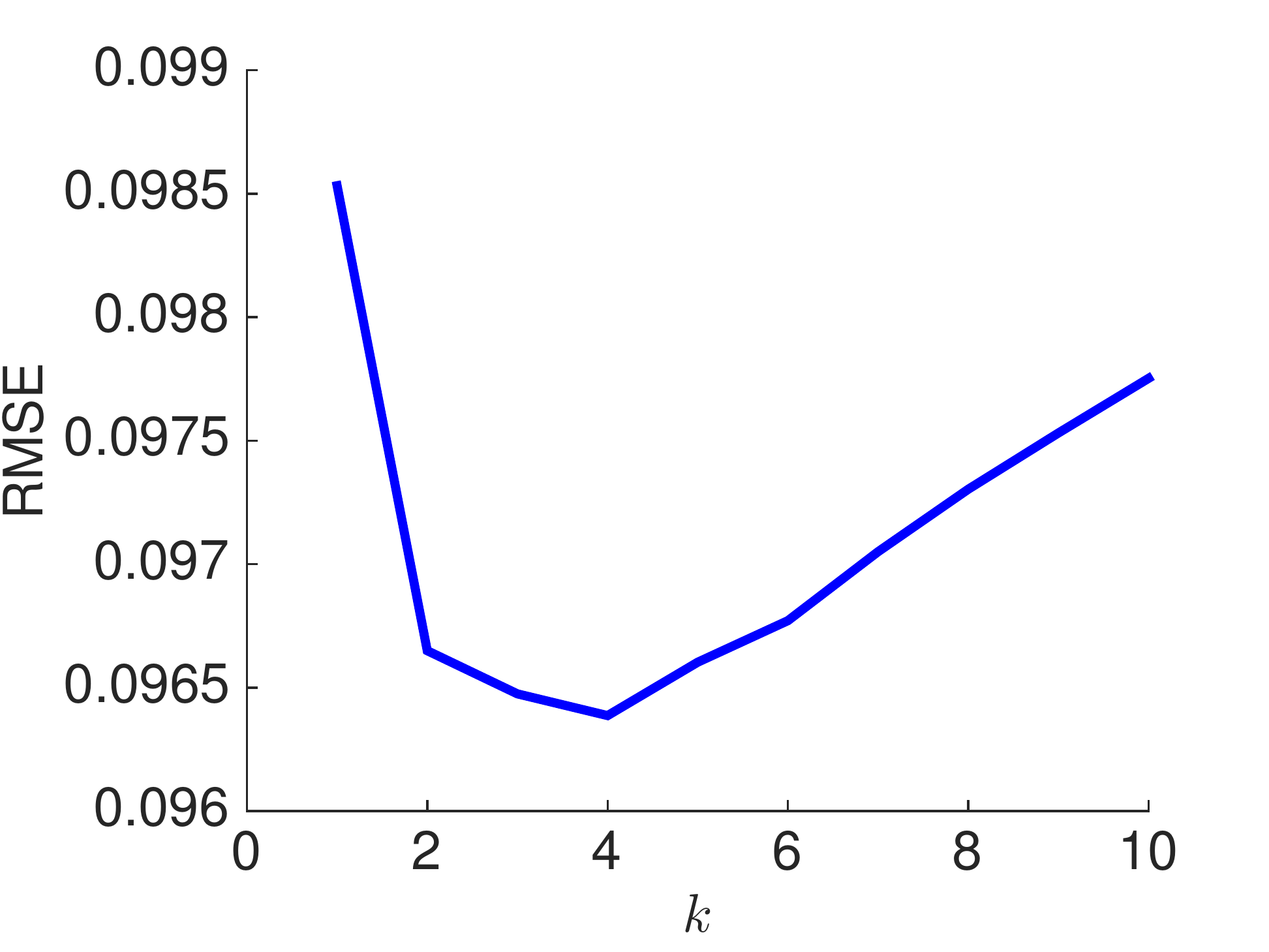}} ~
\subfloat[Iterations vs. $k$]{\includegraphics[width=\figwidth]{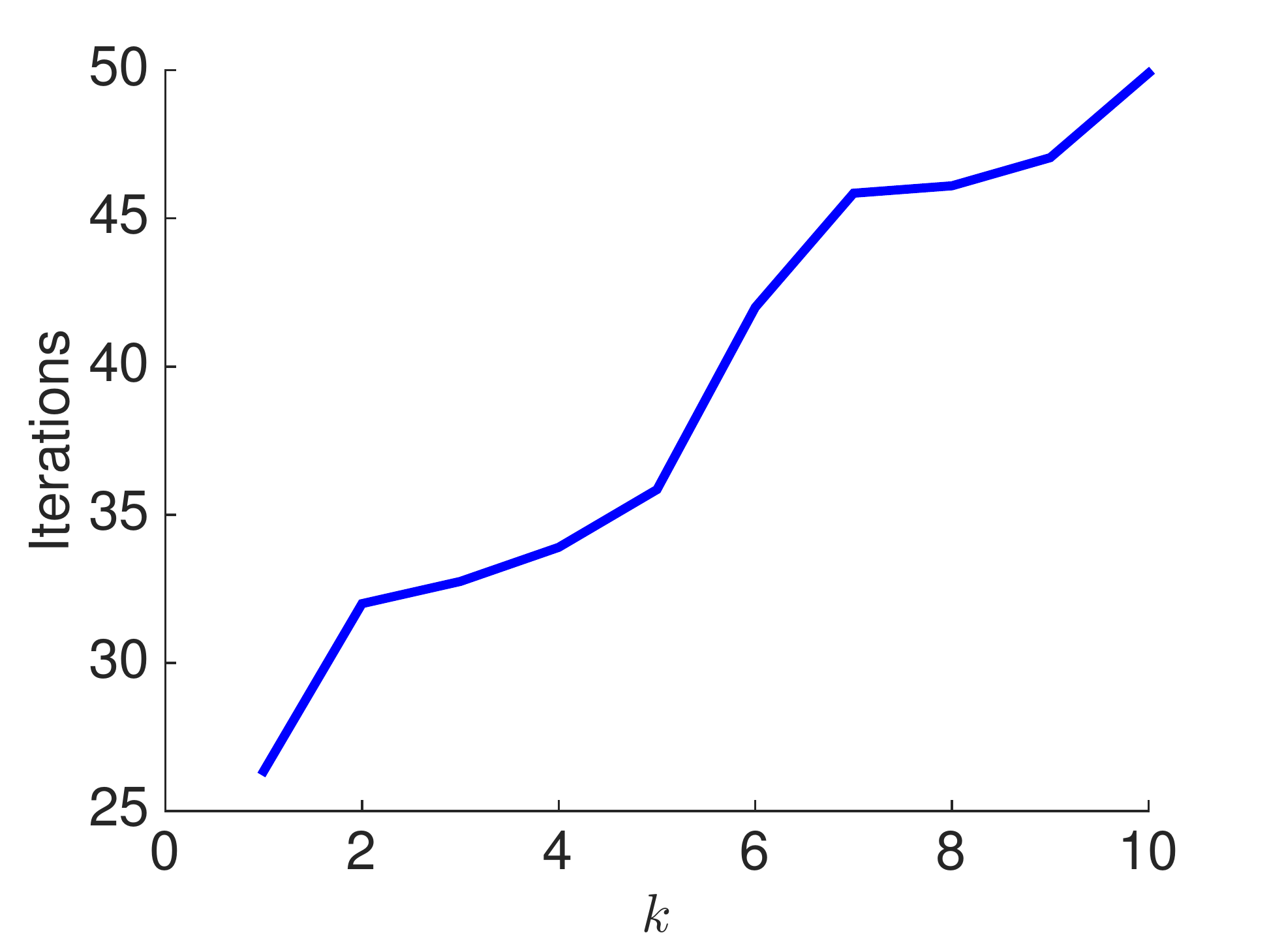}} ~
\subfloat[RMSE vs. Iterations]{\includegraphics[width=\figwidth]{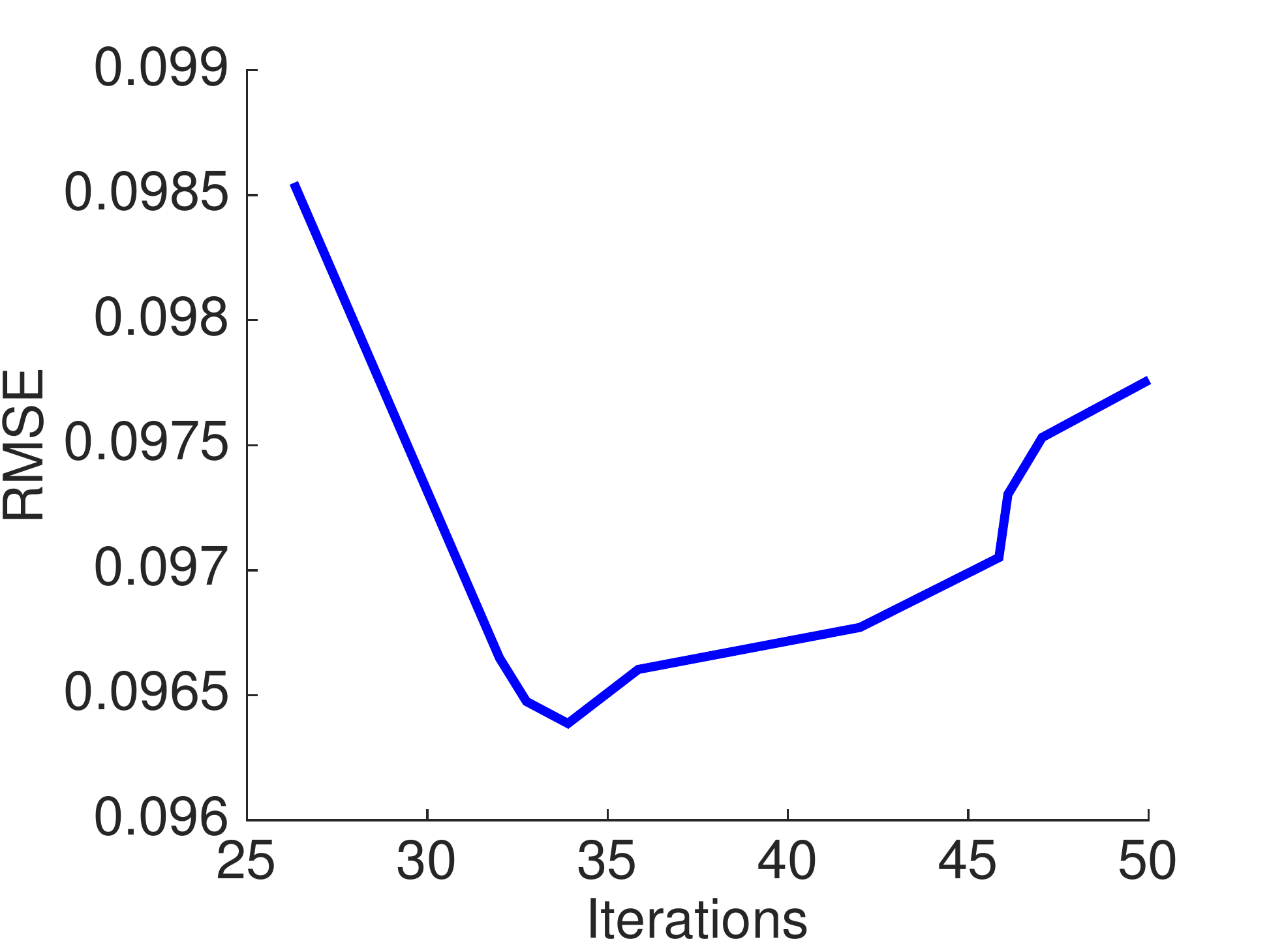}} ~
\caption{\textbf{Computational-statistical tradeoff:} These plots show algorithmic performances for the Couple image with mean 10 corrupted by Poisson noise using the proposed $\log$-$\exp$ \protect\eqref{eq:poisson_final}. A dip in the RMSE vs. $k$ plot (a) and a generally increasing behavior in algorithm Iterations vs. $k$ (b) shows us that there is a dip in RMSE vs. Iterations (c). This plot in (c) suggests how well (in terms of RMSE) we can expect to perform with the proposed reparameterizations when given a computational budget.}
\label{fig:compstat}
\end{center}
\end{figure}

\subsection{Discussion}

In these experiments, the values of $\tau$ were chosen clarivoyantly to minimize the RMSE for each reparameterization considered. In practical settings, we do not have access to ground truth to compute RMSE, but we can use cross-validation to select $\tau$ and the specific $k$ parameter within the family considered to obtain an estimate. If on a time budget, we can leverage the experimental results to eliminate the search over $k$'s:
\begin{itemize}
\item In the Bernoulli case, Table \ref{table:table_bernoulli} in the appendices shows that the value of $k$ that minimizes the risk $\hat{R}(f,\rho)$ for the average of the bank of 10 images is $k=0.3$. 
\item In the Poisson case, Tables \ref{table:table_poisson_mean5} and \ref{table:table_poisson_mean10} in the appendices show that $k=2$ minimizes the average risk $\hat{R}(f,\rho)$ over the 10 images with mean 5, and $k=4$ minimizes the average risk for images with mean 10; to express $k$ as a fraction of the mean we can say $k=0.4 \bar{x^*}$. In practical settings, for a large number of pixels we can roughly estimate $k \approx 0.4 \bar{y}$ to use in the data-fit \eqref{eq:poisson_final}.
\item In the speckle noise case, we have already established that as $k$ grows, the average risk $\hat{R}(f,\rho)$ stabilizes. In practical settings, we choose to use $k = 2$ when applying it to the data-fit \eqref{eq:gamma_final}.
\end{itemize}

\section{Conclusion}
\label{sec:conclusion}
We have introduced a reparameterization framework for three types of noise commonly found in imaging problems, where the reparameterizations proposed are determined by the structure of the negative log-likelihood. This essentially models the image as a parameterized function, allowing us to implicitly modify the regularizer. Image reconstruction is still practical with this change, as we can still use well-studied proximal denoising subroutines in the optimization. Through examples from the exponential family and speckle noise, we see that reparameterization yields gains in RMSE, and has computational implications. Because searching over all possible reparameterization functions is not practical, the functions presented in this paper are by no means ``optimal"; we simply introduce a generalized framework (with example reparameterizations) that may be improved upon with the selection of other functions in future work.

\clearpage
\appendix

\section{Additional numerical results for Bernoulli experiments}

\setcounter{table}{0}
\renewcommand{\thetable}{A\arabic{table}}
\setcounter{figure}{0}
\renewcommand{\thefigure}{A\arabic{figure}}

\begin{table}[ht]
\tiny
\begin{center}
\begin{tabular}{c | c | c c c c c c c c c c}
&&\multicolumn{10}{c}{$k$} \\
Image 		& Baseline & 0.05 & 0.10 & 0.15 & 0.20 & 0.25 & 0.30 & 0.35 & 0.40 & 0.45 & 0.50 \\
\hline
Barbara 		& 0.1766 & 0.1685 & 0.1661 & 0.1643 & 0.1635 & 0.1636 & 0.1632 & \bf{0.1632} & 0.1634 & 0.1635 & 0.1635 \\
Boat 			& 0.1507 & 0.1415 & 0.1391 & 0.1373 & \bf{0.1368} & 0.1375 & 0.1372 & 0.1379 & 0.1376 & 0.1379 & 0.1376 \\
Cameraman 	& 0.1658 & 0.1609 & 0.1585 & \bf{0.1573} & 0.1581 & 0.1590 & 0.1600 & 0.1605 & 0.1610 & 0.1613 & 0.1614 \\
Couple 		& 0.1613 & 0.1531 & 0.1509 & 0.1498 & \bf{0.1494} & 0.1503 & 0.1503 & 0.1508 & 0.1507 & 0.1509 & 0.1509 \\
Fingerprint 	& 0.2207 & 0.2111 & 0.2087 & 0.2056 & 0.2029 & 0.2010 & 0.1998 & 0.1995 & \bf{0.1990} & 0.1991 & 0.1991 \\
Hill 			& 0.1443 & 0.1373 & 0.1351 & 0.1331 & 0.1322 & 0.1316 & 0.1314 & 0.1313 & \bf{0.1313} & 0.1314 & 0.1315 \\
House 		& 0.1148 & 0.1099 & \bf{0.1092} & 0.1104 & 0.1104 & 0.1095 & 0.1102 & 0.1103 & 0.1107 & 0.1111 & 0.1112 \\
Lena 		& 0.1254 & 0.1165 & 0.1153 & 0.1145 & \bf{0.1141} & 0.1159 & 0.1147 & 0.1146 & 0.1151 & 0.1152 & 0.1158 \\
Man 			& 0.1542 & 0.1444 & 0.1437 & 0.1424 & 0.1412 & 0.1412 & 0.1411 & \bf{0.1404} & 0.1406 & 0.1408 & 0.1409 \\
Peppers 		& 0.1754 & 0.1678 & 0.1651 & 0.1638 & 0.1628 & 0.1627 & 0.1630 & 0.1628 & \bf{0.1626} & 0.1629 & 0.1631 \\
\hline
Average		& 0.1589 & 0.1511 & 0.1492 & 0.1479 & 0.1471 & 0.1472 & \bf{0.1471} & 0.1471 & 0.1472 & 0.1474 & 0.1475 \\
\end{tabular}
\end{center}
\caption{\small \textbf{Bernoulli noise data:} RMSE results (all averaged over 20 noise realizations) for each image for the baseline reconstruction \protect\eqref{eq:baseline} and for 10 different reparameterizations from $k = 0.05$ to $0.50$ with the proposed \protect\eqref{eq:bernoulli_final}. Each row has the best performing RMSE in bold. The last row averages the columns, which is the relative risk $\hat{R}(f,\rho)$ defined by \protect\eqref{eq:relative_risk}, and these values are plotted in Figure \ref{fig:bern_plots}. We see that the baseline estimate performs worse than the best reparameterization (in fact, for all the $k$ values chosen, the baseline is overall worse). }
\label{table:table_bernoulli}
\end{table}

\setlength{\figwidth}{0.2\textwidth}
\begin{figure}[ht]
\begin{center}
\subfloat[Truth]{\includegraphics[width=\figwidth]{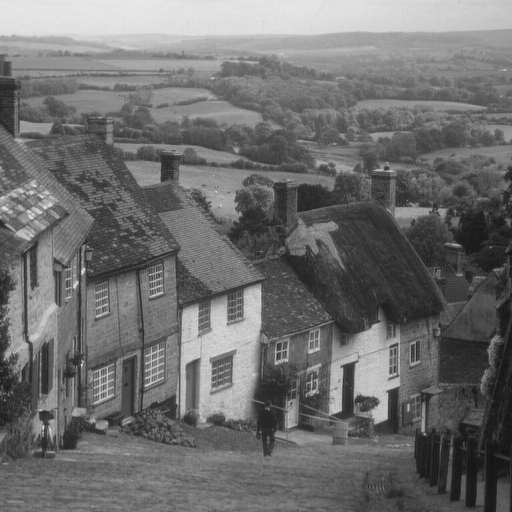}} ~
\subfloat[Observations]{\includegraphics[width=\figwidth]{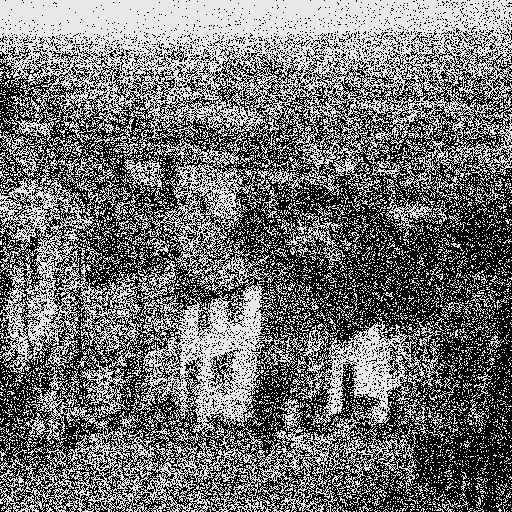}} ~
\subfloat[Baseline recovery]{\includegraphics[width=\figwidth]{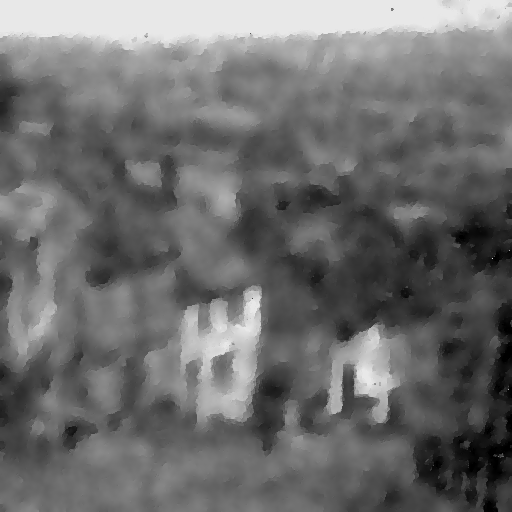}} ~
\subfloat[Reparameterized recovery ($k = 0.40$)]{\includegraphics[width=\figwidth]{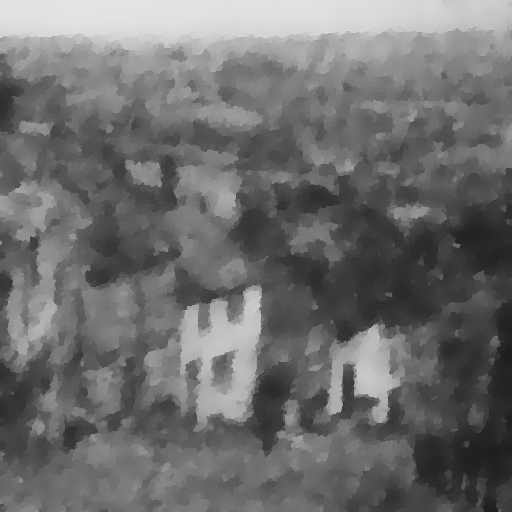}} \\
\subfloat[Truth]{\includegraphics[width=\figwidth]{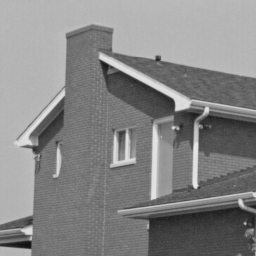}} ~
\subfloat[Observations]{\includegraphics[width=\figwidth]{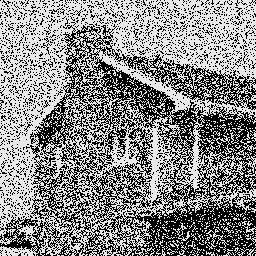}} ~
\subfloat[Baseline recovery]{\includegraphics[width=\figwidth]{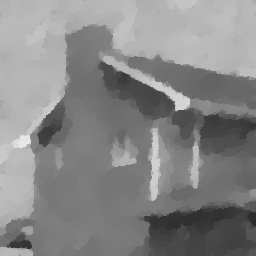}} ~
\subfloat[Reparameterized recovery ($k = 0.10$)]{\includegraphics[width=\figwidth]{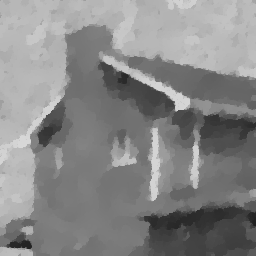}}
\caption{\small \textbf{Bernoulli noise images:} Reconstruction results from the ``Hill" image (top row); reconstruction results from the ``House" image (bottom row). These images are all displayed on the same intensity scale (from 0 to 1.1$\times \max(x^*)$).}
\label{fig:bernoulli_images}
\end{center}
\end{figure}

\clearpage
\section{Additional numerical results for Poisson experiments (mean 5)}

\setcounter{table}{0}
\renewcommand{\thetable}{B\arabic{table}}
\setcounter{figure}{0}
\renewcommand{\thefigure}{B\arabic{figure}}

\begin{table}[ht]
\tiny
\begin{center}
\begin{tabular}{c | c | c c c c c c c c c c}
&&\multicolumn{10}{c}{$k$} \\
Image 		& Baseline & 0.5 & 1.0 & 1.5 & 2.0 & 2.5 & 3.0 & 3.5 & 4.0 & 4.5 & 5.0 \\
\hline
Barbara 		& 0.1472 & 0.1404 & 0.1385 & 0.1375 & 0.1369 & 0.1365 & \bf{0.1365} & 0.1367 & 0.1369 & 0.1369 & 0.1370 \\
Boat 			& 0.1091 & 0.1045 & 0.1030 & \bf{0.1029} & 0.1030 & 0.1035 & 0.1038 & 0.1042 & 0.1043 & 0.1046 & 0.1048 \\
Cameraman 	& 0.1228 & 0.1115 & 0.1107 & \bf{0.1105} & 0.1109 & 0.1114 & 0.1120 & 0.1129 & 0.1134 & 0.1140 & 0.1144 \\
Couple 		& 0.1171 & 0.1130 & 0.1119 & 0.1121 & \bf{0.1116} & 0.1120 & 0.1123 & 0.1126 & 0.1129 & 0.1134 & 0.1136 \\
Fingerprint 	& 0.1698 & 0.1599 & 0.1557 & 0.1528 & 0.1515 & \bf{0.1513} & 0.1516 & 0.1521 & 0.1520 & 0.1525 & 0.1527 \\
Hill 			& 0.1110 & 0.1039 & 0.1019 & 0.1012 & \bf{0.1009} & 0.1010 & 0.1011 & 0.1015 & 0.1016 & 0.1018 & 0.1022 \\
House 		& 0.0931 & 0.0888 & 0.0868 & 0.0854 & 0.0854 & \bf{0.0852} & 0.0856 & 0.0856 & 0.0859 & 0.0860 & 0.0862 \\
Lena 		& 0.0889 & 0.0877 & 0.0850 & 0.0838 & \bf{0.0838} & 0.0840 & 0.0846 & 0.0849 & 0.0852 & 0.0856 & 0.0860 \\
Man 			& 0.1154 & 0.1091 & 0.1066 & 0.1057 & \bf{0.1050} & 0.1050 & 0.1053 & 0.1054 & 0.1054 & 0.1056 & 0.1060 \\
Peppers 		& 0.1131 & 0.1061 & 0.1046 & 0.1044 & \bf{0.1044} & 0.1050 & 0.1055 & 0.1059 & 0.1063 & 0.1068 & 0.1074 \\
\hline
Average 		& 0.1187 & 0.1125 & 0.1105 & 0.1096 & \bf{0.1093} & 0.1095 & 0.1098 & 0.1102 & 0.1104 & 0.1107 & 0.1110 \\
\end{tabular}
\end{center}
\caption{\small \textbf{Poisson noise data (mean 5):} RMSE results (all averaged over 20 noise realizations) for each image for the baseline reconstruction \protect\eqref{eq:baseline} and for 10 different reparameterizations from $k = 0.5$ to $5.0$ with the proposed $\log$-$\exp$ \protect\eqref{eq:poisson_final}. Each row has the best performing RMSE in bold. The last row averages the columns, which is the relative risk $\hat{R}(f,\rho)$ defined by \protect\eqref{eq:relative_risk}, and these values are plotted in Figure \ref{fig:poiss_mean5_plots}. We see that the baseline estimate performs worse than the best reparameterization (in fact, for all the $k$ values chosen, the baseline is overall worse).}
\label{table:table_poisson_mean5}
\end{table}

\setlength{\figwidth}{0.2\textwidth}
\begin{figure}[ht]
\begin{center}
\subfloat[Truth]{\includegraphics[width=\figwidth]{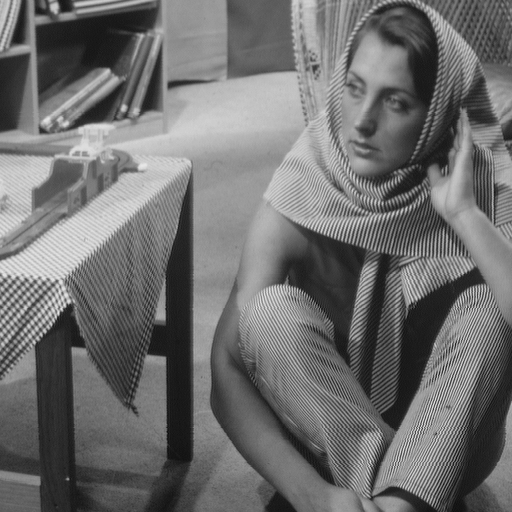}} ~
\subfloat[Observations]{\includegraphics[width=\figwidth]{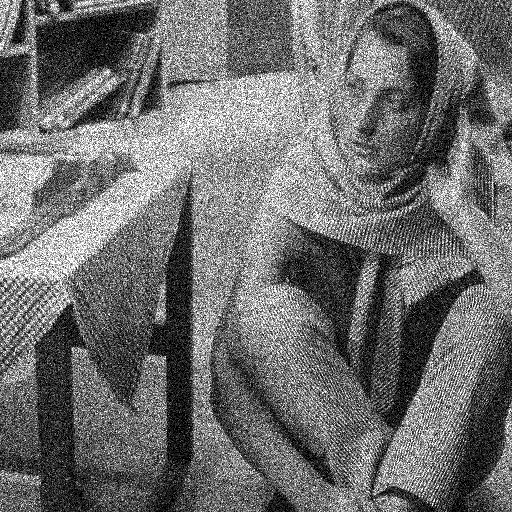}} ~
\subfloat[Baseline recovery]{\includegraphics[width=\figwidth]{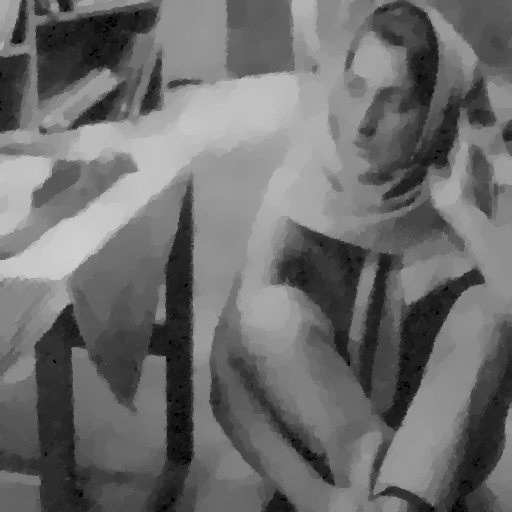}} ~
\subfloat[Reparameterized recovery ($k = 3.0$)]{\includegraphics[width=\figwidth]{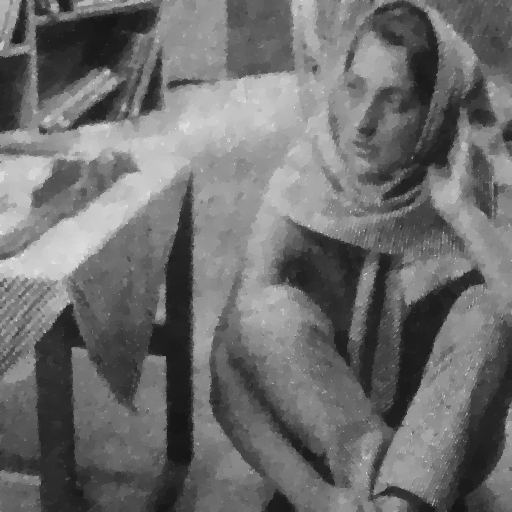}} \\
\subfloat[Truth]{\includegraphics[width=\figwidth]{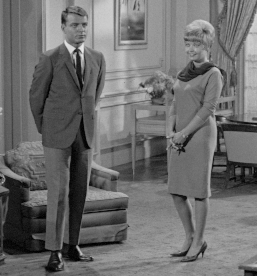}} ~
\subfloat[Observations]{\includegraphics[width=\figwidth]{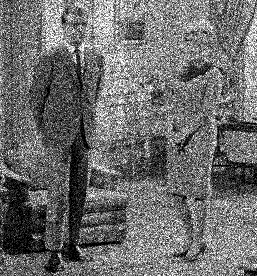}} ~
\subfloat[Baseline recovery]{\includegraphics[width=\figwidth]{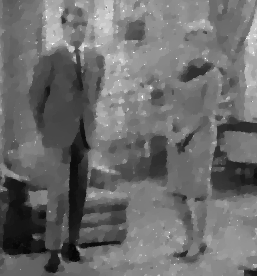}} ~
\subfloat[Reparameterized recovery ($k=2.0$)]{\includegraphics[width=\figwidth]{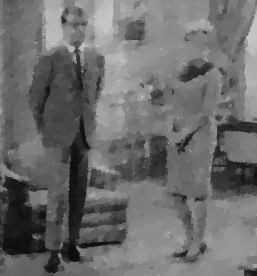}} 
\caption{\small \textbf{Poisson noise images (mean 5):} Reconstruction results from the ``Barbara" image (top row); reconstruction results from the (zoomed-in) ``Couple" image (bottom row). These images are all displayed on the same intensity scale (from 0 to 1.1$\times \max(x^*)$). }
\label{fig:poisson_images_mean5}
\end{center}
\end{figure}

\clearpage
\section{Additional numerical results for Poisson experiments (mean 10)}

\setcounter{table}{0}
\renewcommand{\thetable}{C\arabic{table}}
\setcounter{figure}{0}
\renewcommand{\thefigure}{C\arabic{figure}}

\begin{table}[ht]
\tiny
\begin{center}
\begin{tabular}{c | c | c c c c c c c c c c}
&&\multicolumn{10}{c}{$k$} \\
Image 		& Baseline & 1 & 2 & 3 & 4 & 5 & 6 & 7 & 8 & 9 & 10 \\
\hline
Barbara 		& 0.1369 & 0.1271 & 0.1246 & 0.1228 & 0.1221 & \bf{0.1220} & 0.1220 & 0.1221 & 0.1222 & 0.1223 & 0.1224 \\
Boat 			& 0.0954 & 0.0906 & 0.0890 & 0.0889 & 0.0888 & \bf{0.0888} & 0.0891 & 0.0893 & 0.0895 & 0.0897 & 0.0900 \\
Cameraman 	& 0.1022 & 0.0955 & 0.0924 & \bf{0.0924} & 0.0928 & 0.0933 & 0.0938 & 0.0942 & 0.0947 & 0.0951 & 0.0955 \\
Couple 		& 0.1031 & 0.0985 & 0.0966 & 0.0965 & \bf{0.0964} & 0.0966 & 0.0968 & 0.0971 & 0.0973 & 0.0975 & 0.0978 \\
Fingerprint 	& 0.1347 & 0.1301 & 0.1280 & 0.1263 & \bf{0.1255} & 0.1255 & 0.1257 & 0.1259 & 0.1260 & 0.1262 & 0.1264 \\
Hill 			& 0.0941 & 0.0905 & 0.0895 & 0.0887 & 0.0884 & \bf{0.0882} & 0.0884 & 0.0885 & 0.0887 & 0.0889 & 0.0891 \\
House 		& 0.0745 & 0.0717 & 0.0711 & 0.0709 & 0.0708 & \bf{0.0704} & 0.0705 & 0.0706 & 0.0707 & 0.0708 & 0.0709 \\
Lena 		& 0.0755 & 0.0727 & 0.0717 & \bf{0.0716} & 0.0717 & 0.0718 & 0.0721 & 0.0723 & 0.0726 & 0.0729 & 0.0731 \\
Man 			& 0.0975 & 0.0943 & 0.0925 & 0.0916 & 0.0910 & \bf{0.0908} & 0.0909 & 0.0910 & 0.0911 & 0.0913 & 0.0914 \\
Peppers 		& 0.0937 & 0.0878 & 0.0873 & 0.0865 & \bf{0.0865} & 0.0866 & 0.0869 & 0.0872 & 0.0875 & 0.0878 & 0.0882 \\
\hline
Average 		& 0.1008 & 0.0959 & 0.0943 & 0.0936 & \bf{0.0934} & 0.0934 & 0.0936 & 0.0938 & 0.0940 & 0.0942 & 0.0945 \\
\end{tabular}
\end{center}
\caption{\small \textbf{Poisson noise data (mean 10):} RMSE results (all averaged over 20 noise realizations) for each image for the baseline reconstruction \protect\eqref{eq:baseline} and for 10 different reparameterizations from $k = 1$ to $10$ with the proposed $\log$-$\exp$ \protect\eqref{eq:poisson_final}. Each row has the best performing RMSE in bold. The last row averages the columns, which is the relative risk $\hat{R}(f,\rho)$ defined by \protect\eqref{eq:relative_risk}, and these values are plotted in Figure \ref{fig:poiss_mean10_plots}. We see that the baseline estimate performs worse than the best reparameterization (in fact, for all the $k$ values chosen, the baseline is overall worse).}
\label{table:table_poisson_mean10}
\end{table}

\setlength{\figwidth}{0.2\textwidth}
\begin{figure}[ht]
\begin{center}
\subfloat[Truth]{\includegraphics[width=\figwidth]{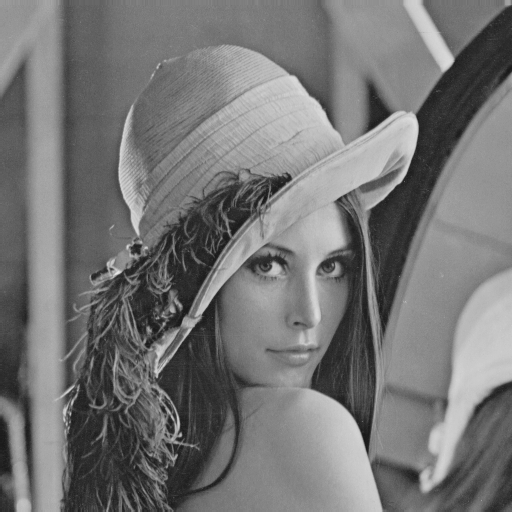}} ~
\subfloat[Observations]{\includegraphics[width=\figwidth]{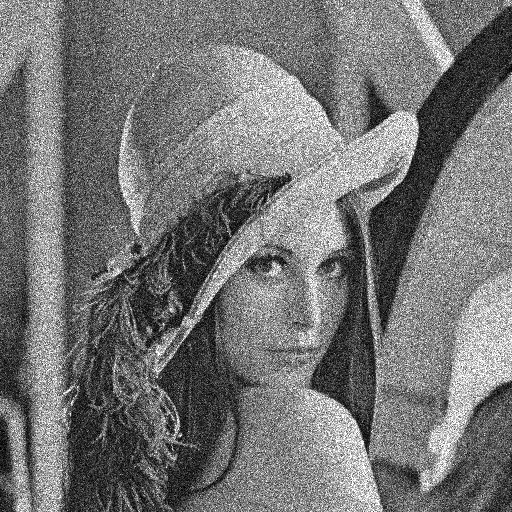}} ~
\subfloat[Baseline recovery]{\includegraphics[width=\figwidth]{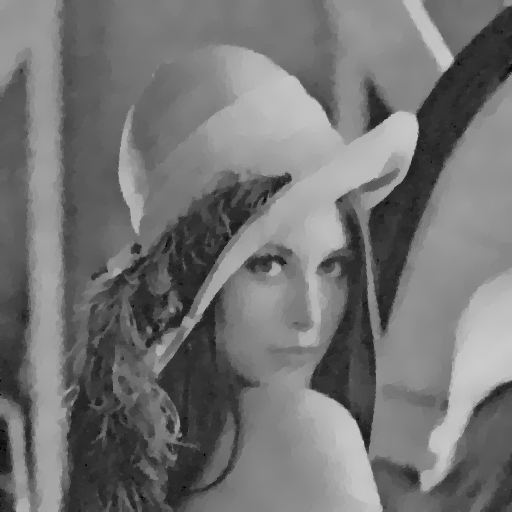}} ~
\subfloat[Reparameterized recovery ($k=3$)]{\includegraphics[width=\figwidth]{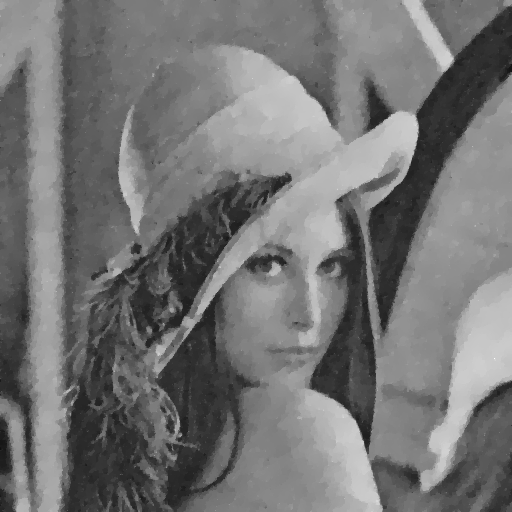}} \\
\subfloat[Truth]{\includegraphics[width=\figwidth]{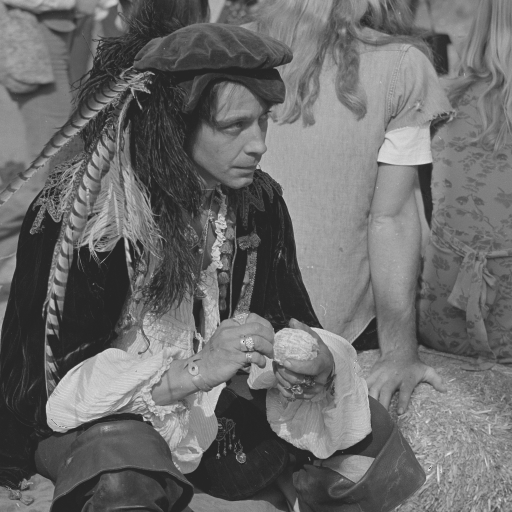}} ~
\subfloat[Observations]{\includegraphics[width=\figwidth]{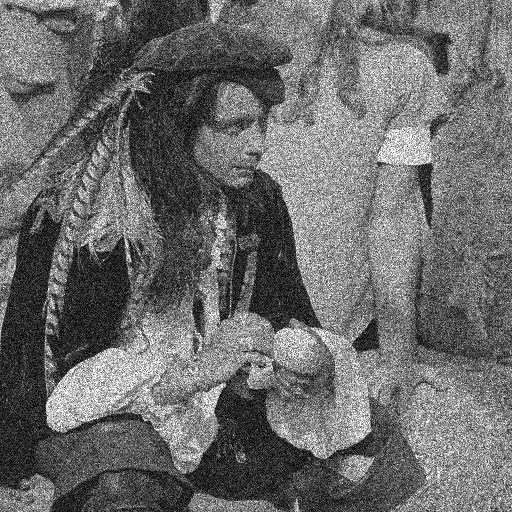}} ~
\subfloat[Baseline recovery]{\includegraphics[width=\figwidth]{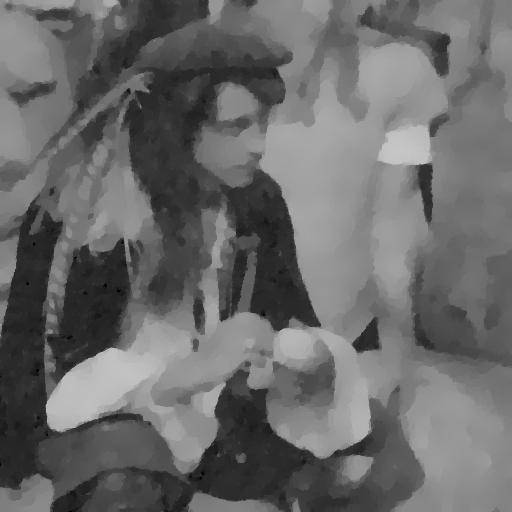}} ~
\subfloat[Reparameterized recovery ($k=5$)]{\includegraphics[width=\figwidth]{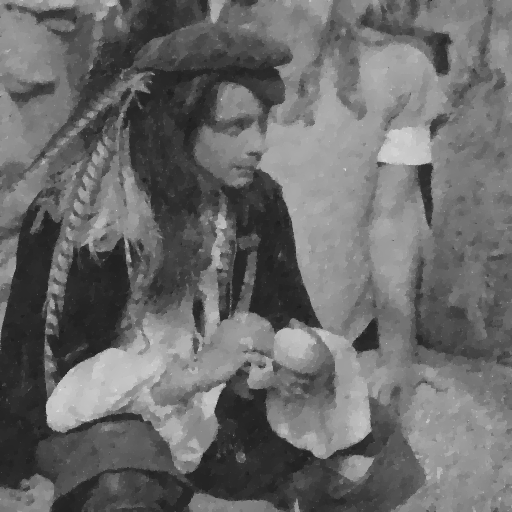}} 
\caption{\small \textbf{Poisson noise images (mean 10):} Reconstruction results from the ``Lena" image (top row); reconstruction results from the ``Man" image (bottom row). These images are all displayed on the same intensity scale (from 0 to 1.1$\times \max(x^*)$).}
\label{fig:poisson_images_mean10}
\end{center}
\end{figure}

\clearpage
\section{Numerical results for speckle noise experiments}

\setcounter{table}{0}
\renewcommand{\thetable}{D\arabic{table}}
\setcounter{figure}{0}
\renewcommand{\thefigure}{D\arabic{figure}}

\begin{table}[ht]
\tiny
\begin{center}
\begin{tabular}{c | c c | c c | c c | c c}
& \multicolumn{2}{c}{Mean 3, $S=3$} & \multicolumn{2}{c}{Mean 3, $S=5$} & \multicolumn{2}{c}{Mean 5, $S=3$} & \multicolumn{2}{c}{Mean 5, $S=5$}\\
Image 		& Baseline & $k = 2$ & Baseline & $k = 2$ & Baseline & $k = 2$ & Baseline & $k = 2$\\
\hline 
Barbara		& 0.1717 & 0.1491 & 0.1664 & 0.1397 & 0.1969 & 0.1493 & 0.1810 & 0.1399 \\
Boat			& 0.1373 & 0.1183 & 0.1221 & 0.1056 & 0.1568 & 0.1184 & 0.1385 & 0.1059 \\
Cameraman	& 0.1530 & 0.1312 & 0.1398 & 0.1155 & 0.1716 & 0.1313 & 0.1594 & 0.1161 \\
Couple		& 0.2176 & 0.1273 & 0.1694 & 0.1140 & 0.1794 & 0.1271 & 0.1656 & 0.1139 \\
Fingerprint	& 0.2261 & 0.1862 & 0.1980 & 0.1607 & 0.2595 & 0.1861 & 0.2120 & 0.1606 \\
Hill			& 0.1580 & 0.1157 & 0.1446 & 0.1047 & 0.1649 & 0.1160 & 0.1517 & 0.1049 \\
House		& 0.1504 & 0.1023 & 0.1362 & 0.0896 & 0.1441 & 0.1012 & 0.1369 & 0.0891 \\
Lena			& 0.1331 & 0.0981 & 0.1219 & 0.0870 & 0.1413 & 0.0982 & 0.1291 & 0.0868 \\
Man			& 0.1530 & 0.1218 & 0.1412 & 0.1092 & 0.1656 & 0.1218 & 0.1543 & 0.1093 \\
Peppers		& 0.1489 & 0.1256 & 0.1329 & 0.1084 & 0.1764 & 0.1247 & 0.1629 & 0.1077 \\
\hline
Average		& 0.1649 & 0.1276 & 0.1472 & 0.1134 & 0.1757 & 0.1274 & 0.1591 & 0.1134 \\
\end{tabular}
\end{center}
\caption{\small \textbf{Speckle noise data:} RMSE results (all averaged over 20 noise realizations) for each image for the baseline reconstruction \protect\eqref{eq:baseline} and for the reparameterization with $k = 2$ with the proposed \protect\eqref{eq:gamma_final}; we only state the results for $k = 2$ because Figure \ref{fig:gamma_plots} suggests that the RMSE levels off as $k$ increases (and for the sake of space). The last row averages the columns, which is the relative risk $\hat{R}(f,\rho)$ defined by \protect\eqref{eq:relative_risk}, and these values are plotted in Figure \ref{fig:gamma_plots}. We see that the baseline estimate performs worse than the best reparameterization at the leveled-off valued of $k = 2$.}
\label{table:table_gamma}
\end{table}

\setlength{\figwidth}{0.2\textwidth}
\begin{figure}[ht]
\begin{center}
\subfloat[Truth]{\includegraphics[width=\figwidth]{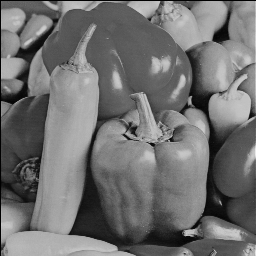}} ~
\subfloat[Observations]{\includegraphics[width=\figwidth]{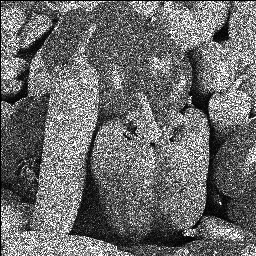}} ~
\subfloat[Baseline recovery]{\includegraphics[width=\figwidth]{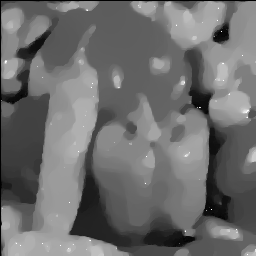}} ~
\subfloat[Reparameterized recovery ($k=2$)]{\includegraphics[width=\figwidth]{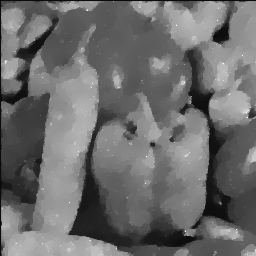}} \\
\subfloat[Truth]{\includegraphics[width=\figwidth]{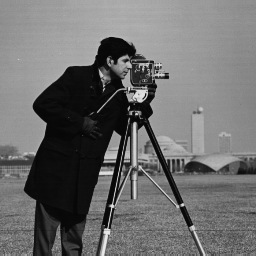}} ~
\subfloat[Observations]{\includegraphics[width=\figwidth]{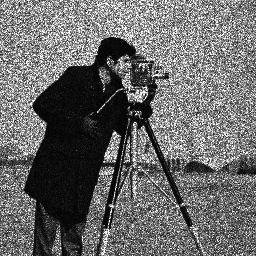}} ~
\subfloat[Baseline recovery]{\includegraphics[width=\figwidth]{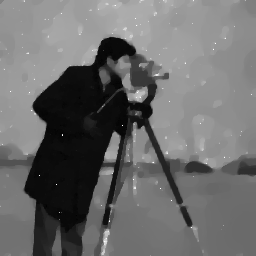}} ~
\subfloat[Reparameterized recovery ($k=2$)]{\includegraphics[width=\figwidth]{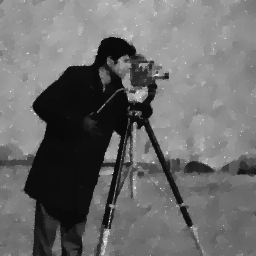}}
\caption{\small \textbf{Speckle noise images (mean 3):} Reconstruction results from the ``Peppers" image with $S=3$ (top row); reconstruction results from the ``Cameraman" image with $S=5$ (bottom row). These images are all displayed on the same intensity scale (from 0 to 1.1$\times \max(x^*)$).}
\label{fig:gamma_images_mean3}
\end{center}
\end{figure}

\setlength{\figwidth}{0.2\textwidth}
\begin{figure}[ht]
\begin{center}
\subfloat[Truth]{\includegraphics[width=\figwidth]{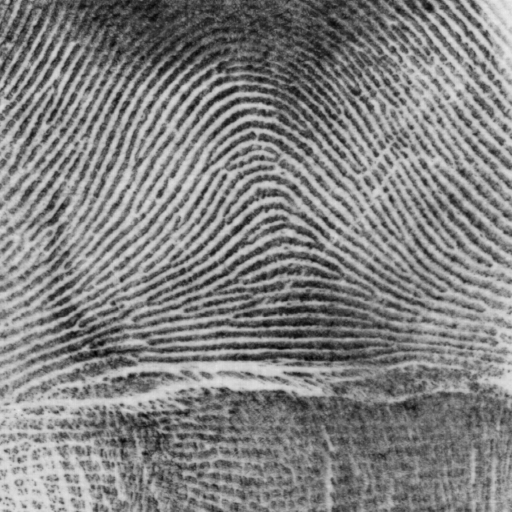}} ~
\subfloat[Observations]{\includegraphics[width=\figwidth]{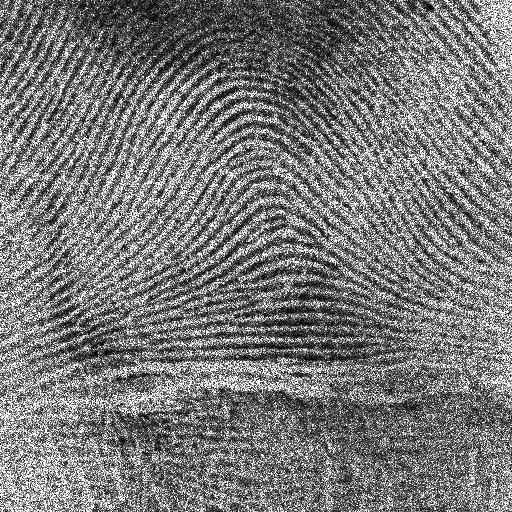}} ~
\subfloat[Baseline recovery]{\includegraphics[width=\figwidth]{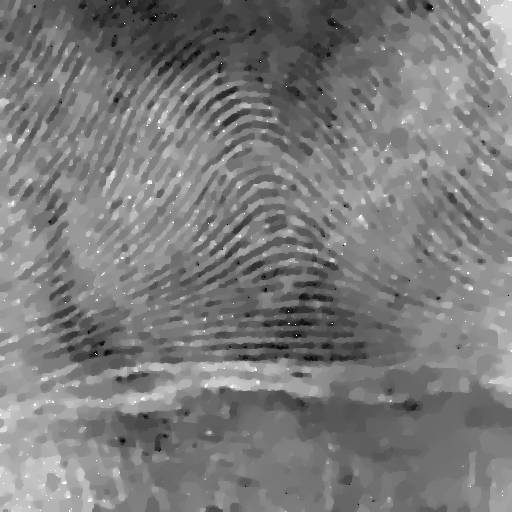}} ~
\subfloat[Reparameterized recovery ($k=2$)]{\includegraphics[width=\figwidth]{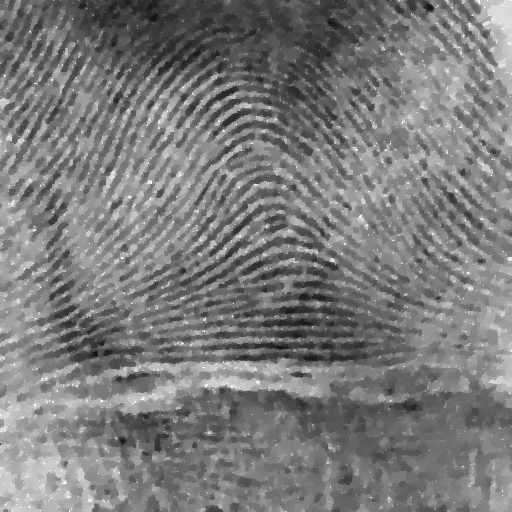}} \\
\subfloat[Truth]{\includegraphics[width=\figwidth]{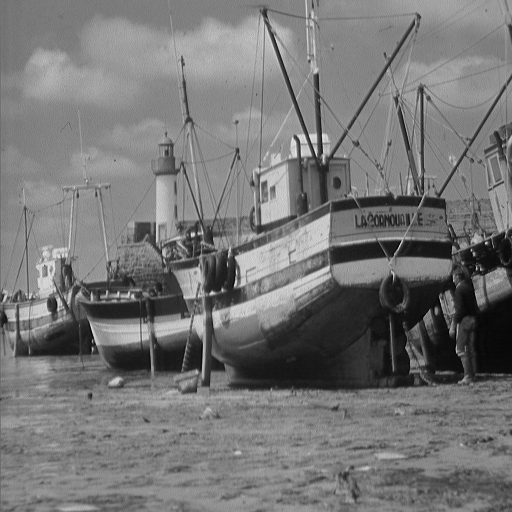}} ~
\subfloat[Observations]{\includegraphics[width=\figwidth]{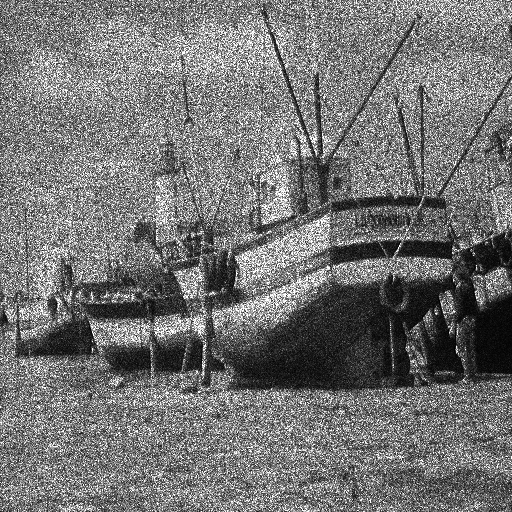}} ~
\subfloat[Baseline recovery]{\includegraphics[width=\figwidth]{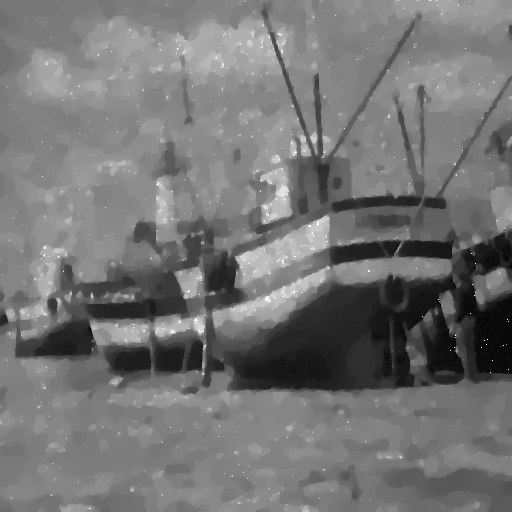}} ~
\subfloat[Reparameterized recovery ($k=2$)]{\includegraphics[width=\figwidth]{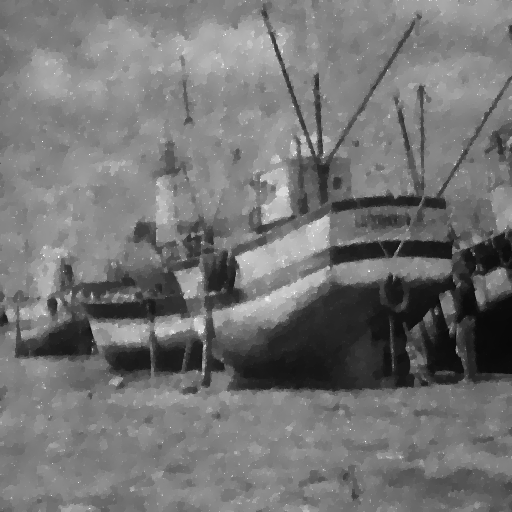}}
\caption{\small \textbf{Speckle noise images (mean 5):} Reconstruction results from the ``Fingerprint" image with $S=3$ (top row); reconstruction results from the ``Boat" image with $S=5$ (bottom row). These images are all displayed on the same intensity scale (from 0 to 1.1$\times \max(x^*)$).}
\label{fig:gamma_images_mean5}
\end{center}
\end{figure}

\newpage
\bibliographystyle{siam}
\bibliography{refs}

\end{document}